\documentclass[smallextended,envcountsect]{svjour3}
\smartqed  
\usepackage{diagbox}
\usepackage{amsmath,multirow,booktabs,subfigure}
\usepackage{amssymb}
\usepackage[mathscr]{eucal}
\usepackage{graphicx}
\usepackage{color}
\usepackage{empheq}
\usepackage{graphics}
\usepackage{epstopdf}
\usepackage{array}
\usepackage{makecell}
\usepackage[justification=centering]{caption}
\usepackage[latin1]{inputenc}
\numberwithin{equation}{section}
\newtheorem{assumption}{\sc\bf Assumption}[section]
\newtheorem{algorithm}{\sc\bf Algorithm}[section]

\begin{document}

\title{ Tseng Splitting Method with Double Inertial Steps for Solving Monotone Inclusion Problems}

\author{ Zhong-bao Wang\and Zhen-yin Lei\and Xin Long\\ \and Zhang-you Chen}

\institute{Zhong-bao Wang,  Corresponding author \at
1. Department of Mathematics, Southwest Jiaotong University\\
 Chengdu, Sichuan 611756,  China\\
  2. National Engineering Laboratory of Integrated
Transportation Big Data Application Technology\\ Chengdu, Sichuan 611756, China\\
3. School of Mathematical Sciences, University of Electronic Science and Technology of China\\
             Chengdu, Sichuan, 611731, China\\
              zhongbaowang@hotmail.com
               \and Zhen-yin Lei \at
              1. Department of Mathematics, Southwest Jiaotong University\\
              Chengdu, Sichuan 611756,  China\\
              2. National Engineering Laboratory of Integrated
              Transportation Big Data Application Technology\\ Chengdu, Sichuan 611756, China\\
              lzy0615edu@163.com
              \and Xin Long  \at
           1. Department of Mathematics, Southwest Jiaotong University\\
 Chengdu, Sichuan 611756,  China\\
  2. National Engineering Laboratory of Integrated
Transportation Big Data Application Technology\\ Chengdu, Sichuan 611756, China\\
              1065123400@qq.com
              \and
              Zhang-you Chen \at
 1. Department of Mathematics, Southwest Jiaotong University\\
 Chengdu, Sichuan 611756,  China\\
 2. National Engineering Laboratory of Integrated
Transportation Big Data Application Technology\\ Chengdu, Sichuan 611756, China\\
              zhangyouchen@swjtu.edu.cn}

\date{Received: date / Accepted: date}
\maketitle
\noindent{\bf Abstract}
\par In this paper, based on double inertial extrapolation steps strategy and relaxed techniques, we introduce a new Tseng splitting method with double inertial  extrapolation steps and self-adaptive step sizes for solving monotone inclusion problems in real Hilbert spaces. Under mild and standard assumptions, we establish successively the weak convergence, nonasymptotic $O(\frac{1}{\sqrt{n}})$ convergence rate, strong convergence and linear convergence rate of the proposed algorithm. Finally, several numerical experiments are provided to illustrate the performance and theoretical outcomes of our algorithm. \\


\section{Introduction}
\setcounter{equation}{0}\qquad

Let $H$ be a real Hilbert space with the inner product $\left\langle { \cdot , \cdot } \right\rangle $ and the induced norm $\left\|  \cdot  \right\|$. The monotone inclusion problem (MIP) is as follows
\begin{equation}
	\label{1.1}
    {\rm find}~~x^* \in H~~{\rm such~that}~~0 \in (A + B)x^*.
\end{equation}
where  $A:H \to H$ is a single mapping and $B:H \to {2^H}$ is a multivalued mapping. The solution set is denoted by $\Omega := {(A + B)^{ - 1}}(0)$.
 \par The monotone inclusion problem has drawn much attention because it provides a broad unifying frame for variational inequalities, convex minimization problems, split feasibility problems and equilibrium problems, and has been applied to solve several real-world problems from machine learning, signal processing and image restoration, see \cite{AB,BH,CD,CW,DS,RF,TYCR}.
\par One of famous methods for solving MIP (\ref{1.1}) is forward-backward splitting method, which was introduced by Passty \cite{PG} and Lions et al. \cite{LP}. This method generates an iterative sequence $\{x_n\}$ in following way
\begin{equation}
	\label{1.2}
	{x_{n + 1}} = {(I + {\lambda _n}B)^{ - 1}}(I - {\lambda _n}A){x_n}.
\end{equation}
where the mapping $A$ is $\frac{1}{L}$-co-coercive, $B$ is maximal monotone, $I$ is an identity mapping on $H$ and $\lambda_n>0$. The operator $(I - {\lambda _n}A)$  is called an forward operator and ${(I + {\lambda _n}B)^{ - 1}}$ is said to be a backward operator.
\par Tseng \cite{TP} proposed a modified forward-backward splitting method (also known as Tseng splitting algorithm), whose iterative formula is as follows
\begin{equation}
	\label{1.3}
	\left\{ \begin{array}{l}
		{y_n} = {(I + {\lambda _n}B)^{ - 1}}(I - {\lambda _n}A){x_n}\\
		{x_{n + 1}} = {y_n} - {\lambda _n}(A{y_n} - A{x_n}),
	\end{array} \right.
\end{equation}
where $A$ is $L$-lipschitz continuous and $\{\lambda_n\}\subset (0,1/L)$.  However, the Lipschitz constant of an operator is often unknown or difficult to
estimate in nonlinear problems. To overcome this drawback, Cholamjiak et al. \cite{CV} introduce a relaxed forward-backward splitting method, which uses a
simple step-size rule without the prior knowledge of Lipschitz constant of the operator, for solving MIP \eqref{1.1} and prove the linear convergence rate of the proposed algorithm.
\par In recent year, the inertial method was introduced in \cite{AF}, which can be regarded as a procedure of speeding up the convergence rate of algorithms. Many researchers utilize inertial methods to design algorithm for solving monotone inclusion problems and variational inequalities, see, for example, \cite{AB,CD,CV,CH,CHG,DQ,LN,PTK,SY,TD,YI}. To enhance the numerical efficiency, \c{C}opur et al. \cite{CHG} introduce firstly the double inertial extrapolation steps for solving quasi-variational inequalities in real Hilbert spaces.  Combining relaxation techniques with the inertial methods, Cholamjiak et al. \cite{CV} modify Tseng splitting method to solve MIP \eqref{1.1} in real Hilbert spaces. Very recently, incorporating double inertial extrapolation steps and relaxation techniques, Yao et al. \cite{YI} present a novel subgradients extragradient method to solve variational inequalities,
and prove its strong convergence, weak convergence and linear convergence, respectively. However, the linear convergence of \cite{YI} is obtained under a single inertia rather than double inertias.
\par This paper devotes to further modifying Tseng splitting method for solving MIP \eqref{1.1} in real Hilbert spaces. We obtain successively the weak convergence, nonasymptotic $O(\frac{1}{\sqrt{n}})$ convergence rate, strong convergence and linear convergence rate of the proposed algorithm. Our results obtained in this paper improve the corresponding results in \cite{AB,CV,CH,VA,YI} as follows:
\par$\bullet$ Combining double inertial extrapolation steps strategy and relaxed techniques, we propose a new Tseng splitting method, which include the corresponding methods considered in \cite{AB,CV,YI} as special cases. The two inertial factors in our algorithm are variable sequences, different from the constant inertial factor in \cite{AB,CV,CH,YI}. Especially, when our algorithm is applied to solving variational inequalities, its some parameters have larger choosing interval than the ones of \cite{YI}. In addition, one of our inertial factors can be equal to 1, which is not allowed in the single inertial methods \cite{CV,CH}, which require that the inertial factor must be strictly less than 1.
 \par$\bullet$ We prove the strong convergence, nonasymptotic $O(\frac{1}{\sqrt{n}})$ convergence rate and linear convergence rate of the proposed algorithm.  Note that that the strong convergence does not require to know the modulus of strong monotonicity and the Lipschitz constant in advance. As far as we know, there is no convergence rate results in the literature for methods with double inertial extrapolation steps for solving MIP \eqref{1.1} in infinite-dimensional Hilbert spaces.
  \par$\bullet$ Our algorithm use double inertial extrapolation steps to accelerate the speed of the algorithm. The step sizes of our algorithm are updated by a simple calculation without knowing the Lipschitz constant of the underlying operator. Some numerical experiments show that our algorithm has better efficiency than the corresponding algorithms in \cite{AB,CV,CH,VA,YI}.

The structure of this article is as follows. In section 2, we recall some essential definitions and results which is relate to this paper. In section 3, we present our algorithm and analyze its weak convergence. In section 4, we establish the strong convergence and the linear convergence rate of our method. In section 5, we present some numerical experiments to demonstrate the performance of our algorithm. We give some concluding remarks in section 6.

\section{Preliminaries}
\setcounter{equation}{0}
In this section, we first give some definitions and results that will be used in this paper.
The weak convergence and strong convergence of sequences are denoted by $\rightharpoonup$ and $\to$, respectively.
\begin{definition}
	\label{d2.1}
	The mapping $A:H\rightarrow H$ is called
	\begin{itemize}
\item[(i)] pseudomonotone on $H$ if $\langle Ax, y - x \rangle \geq 0$ implies that $\langle Ay, y - x \rangle \geq 0,~ \forall x,y \in H$;
		\item[(ii)] monotone on $H$ if
		\begin{equation*}
			\langle Ax-Ay, x-y\rangle\geq 0 ,~\forall~ x,y\in H;
		\end{equation*}
  	    \item[(iii)] $\mu$-strongly monotone on $H$ if there exists a positive constant $\mu>0$
	    \begin{equation*}
		\langle Ax-Ay, x-y\rangle\geq \mu {\left\| {x - y} \right\|^2} ,~\forall~ x,y\in H;
	    \end{equation*}
		\item[(iv)]  $L$-lipschitz continuous on $H$ if there exists a scalar $L>0$ satisfying
		\begin{equation*}
			\|Ax-Ay\|\leq L\|x-y\|, ~\forall~ x,y\in H;
		\end{equation*}
\item[(v)] $r$-strongly pseudomonotone on $H$ if there exists a positive constant $r>0$ such that
	    \begin{center}
		 $\langle Ay, x-y\rangle\geq 0$ implies that $\langle Ax, x-y\rangle\geq r {\left\| {x - y} \right\|^2} ,~\forall~ x,y\in H.$
	    \end{center}
	\end{itemize}
\end{definition}

\begin{definition}
	\label{d2.2}
	The graph of $A$ is the set in $H \times H$ defined by
	\begin{equation*}
		Graph(A):=\{ (x,u):x \in H,u \in Ax\} .
	\end{equation*}
Let $C \subset H$ be a nonempty, closed and convex set. The normal cone $N_C(x)$ of $C$ at $x$ is represented by
	\begin{equation*}
	{N_C}(x): = \left\{ \begin{array}{l}
		\{ z \in H:\left\langle {z,y - x} \right\rangle  \le 0,\forall y \in C\} ,~if~x \in C,\\
		\emptyset ,~~~~~~~~~~~~~~~~~~~~~~~~~~~~~~~~~~~~~~~~~otherwise.
	\end{array} \right.
\end{equation*}
The projection of $x \in H$ onto $C$, denoted by $P_C(x)$, is defined as
	\begin{equation*}
	{P_C}(x): = \arg {\min _{y \in C}}\left\| {x - y} \right\|.
\end{equation*}
and is has $\left\langle {x - {P_C}(x),y - {P_C}(x)} \right\rangle  \le 0,\forall y \in C.$\\
The sequence $\{u_n\}$ is $Q$-linear convergence if there is $q \in (0,1)$ such that $\left\| {{u^{k + 1}} - u} \right\| \leq q\left\| {{u^k} - u} \right\|$ for all $k$ large enough.
\end{definition}

\begin{definition}
	\label{d2.3}
	The set-valued mapping $A:H \to {2^H}$ is called
	\begin{itemize}
		\item[(i)] monotone on $H$ if for all $x,y \in H$, $u \in Ax$ and $v \in Ay$ implies that
		\begin{equation*}
			\langle u-v, x-y\rangle\geq 0;
		\end{equation*}
		\item[(ii)] maximal monotone on $H$ if it is monotone and for any
		$(x,u) \in H \times H,\left\langle {u - v,x - y} \right\rangle  \ge 0$ for every $(y,v) \in Graph(A)$ implies that $u \in Ax$;
\item[(iii)]$\mu$-strongly monotone on $H$ if for all $x,y \in H$, $u \in Ax$ and $v \in Ay$ implies that
\begin{equation*}
	\langle u-v, x-y\rangle\geq \mu {\left\| {x - y} \right\|^2} ,~\forall~ x,y\in H.
\end{equation*}
	\end{itemize}
\end{definition}

\begin{lemma}
\label{l2.1}(\cite{MP})Let $\{\varphi_n\},\{\delta_n\}$ and $\{ {\alpha_n}\} $ be sequences in $[0, + \infty )$ such that
\begin{equation*}
	\varphi_{n+1} \le \varphi_n + {\alpha _n}(\varphi_n - \varphi_{n-1}) + \delta_n,\forall n \ge 1,\sum\limits_{n = 1}^\infty  {\delta_n}  <  + \infty
\end{equation*}
and there exists a real number $\alpha$ with $0 \le {\alpha _n} \le \alpha  < 1$ for all $n \in N $.Then the following hold :
\begin{itemize}
  \item[(i)] $\sum {_{n = 1}^{ + \infty }} {[\varphi_n - \varphi_{n-1}]_ + } <  + \infty $ where ${[t]_ + }: = \max \{ t,0\} $;
  \item[(ii)] there exists ${\varphi^*} \in [0, + \infty )$ such that $\mathop {\lim }\limits_{n \to \infty } \varphi_n = {\varphi^*}$.
\end{itemize}
\end{lemma}
\begin{lemma}
	\label{l2.2}(\cite{OZ})Let $C$ be a nonempty set of $H$ and $\{ {x_n}\} $ be a sequence in $H$ such that the following two conditions hold:
	\begin{itemize}
		\item[(i)] for every $x \in C,\mathop {\lim }\limits_{n \to \infty } \left\| {{x_n} - x} \right\|$ exists;
		\item[(ii)] every sequential weak cluster point of $\{ {x_n}\} $ is in $  C$.
	\end{itemize}
Then  $\{ {x_n}\} $ converges weakly to a point in $ C$.
\end{lemma}
\begin{lemma}
	\label{l2.3}\cite{TX} Let $A:H \to H$ be a maximal monotone mapping and $B:H \to {2^H}$ be a Lipschitz continuous and monotone mapping. Then the mapping $A + B$ is a maximal monotone mapping.
\end{lemma}
\begin{lemma}
	\label{l2.4}\cite{LQ} Let $\{a_n\}$ and $\{b_n\}$ be nonnegative real numbers sequences for which there exists $0 \le q < 1$, so that
\begin{equation*}
		a_{n + 1} \leq q{a_n} + b_n~~for~every~n \in N.
\end{equation*}
	If $\mathop {\lim }\limits_{n \to \infty } {b_n} = 0$, then $\mathop {\lim }\limits_{n \to \infty } {a_n} = 0$.
\end{lemma}
\begin{lemma}
	\label{l2.5}\cite{XH} Let $\{ {\alpha _n}\}$, $\{a_n\}$, $\{b_n\}$ and $\{c_n\}$ be nonnegative real numbers sequences and there exists $n_0 \in N$ such that
	\begin{equation*}
		{a_{n + 1}} \le (1 - {\alpha _n}){a_n} + {\alpha _n}{b_n} + {c_n}~~\forall n \ge 1,
	\end{equation*}
		where $\{ {\alpha _n}\}$, $\{b_n\}$ and $\{c_n\}$ satisfy the following conditions
		\begin{itemize}
			\item[(i)] $\{\alpha _n\}  \subset (0,1)$ and $\sum\limits_{n = 0}^\infty  {{\alpha _n}}  = \infty ;$
			\item[(ii)]  $\lim {\sup _{n \to \infty }}{b_n} \le 0;$
			\item[(iii)]	${c_n} \ge 0,~\forall n \ge 0~,\sum\limits_{n = 1}^\infty  {{c_n}}  < \infty$.
		\end{itemize}
		Then $\mathop {\lim }\limits_{n \to \infty } {a_n} = 0$.

\end{lemma}
\begin{lemma}
	\label{l2.6}\cite{GT} Let $B:H \to {2^H}$ be a set-valued maximal monotone mapping and $A:H \to H$ is a mapping. Define ${T_\lambda }: = {(I + \lambda B)^{ - 1}}(I - \lambda A),\lambda  > 0$. Then $Fix({T_\lambda }) = {(A + B)^{ - 1}}0$.	
\end{lemma}

\begin{lemma}
	\label{l2.7}(\cite{BH},Corollaty 2.14) For all $x,y\in H$ and $\alpha \in R$, the following equality holds:
	\[{\left\| {\alpha x + (1 - \alpha )y} \right\|^2} = \alpha {\left\| x \right\|^2} + (1 - \alpha ){\left\| y \right\|^2} - \alpha (1 - \alpha ){\left\| {x - y} \right\|^2}.\]	
\end{lemma}
\section{Weak convergence }

\setcounter{equation}{0}
In this section, we introduce the Tseng splitting method with double inertial steps to solve MIP \eqref{1.1} and discuss convergence and convergence rate of the new algorithm. We firstly give the following conditions.
	\begin{itemize}
		\item[($C_1$)] The solution set of the inclusion problem \eqref{1.1} is nonempty,  that is, $\Omega\ne \emptyset $.
		\item[($C_2$)] The mappings $A:H \to H$ is $L$-Lipschitz continuous and monotone and the set-valued mapping $B:H \to {2^H}$ is maximal monotone.
      \item[($C_3$)] The  real sequences $\{\alpha_{n}\}$, $\{\beta_n\}$, $\{\beta_n\}$, $\{a_n\}$, $\{p_n\}$ and $\{\mu_n\}$ satisfy the following conditions
		\item[(i)] $0\le\alpha_{n}\le1$;
		\item[(ii)] $0\le\beta_n\le\beta_{n+1}\le\beta<\frac{{3 + 2\varepsilon-\sqrt {8\varepsilon  + 17} }}{{2\varepsilon }},\varepsilon  \in (1, + \infty )$;
		\item[(iii)] $0<\theta<\theta_n\le\theta_{n+1}\le \dfrac{1}{1+\varepsilon},\varepsilon  \in (1, + \infty )$;
		\item[(iv)] ${a_n} = (1 - {\theta _n}){\beta _n} + {\theta _n}{\alpha _n}$  is a non-decreasing sequence;
		\item[(v)] $\sum\limits_{n = 1}^\infty  {{p_n}}  < \infty $ and $\mathop {\lim }\limits_{n \to \infty } {\mu _n} = 0$.
	\end{itemize}
\begin{remark}
\label{r3.5}
Note that if $\{\alpha_{n}\}$ is a non-decreasing sequence, $\tilde{\delta}\geq 0$ and $\beta_n=\tilde{\delta}\leq \alpha_{1}$, then the condition ($C_3$)(iv) holds naturally. In addition, if we choose $\alpha_n=\frac{1}{5} + \frac{1}{{6 + n}}$, $\beta_n=\frac{1}{6} - \frac{1}{{6 + n}}$ and $\theta_n=\frac{1}{4} - \frac{1}{{6 + n}}$, then the condition ($C_3$)(iv) is also true.	
\end{remark}	

\begin{algorithm}
		\label{alg3.1}
Choose ${x_0},{x_1} \in H$ , $\mu  \in (0,1)$ and $\lambda _1 >0$
\begin{itemize}
\item[Step 1.]Compute
\begin{equation*}
	 \begin{array}{l}
		{w_n} = {x_n} + {\alpha _n}({x_n} - {x_{n - 1}})\\
		{z_n} = {x_n} + {\beta _n}({x_n} - {x_{n - 1}})\\
		{y_n} = {(I + {\lambda _n}B)}^{-1}(I - {\lambda _n}A){w_n}
	\end{array}
\end{equation*}
where
\begin{equation*}
{\lambda _{n + 1}} = \left\{ \begin{array}{l}
	\min \{ \frac{{(\mu_n+\mu) \left\| {{w_n} - {y_n}} \right\|}}{{\left\| {A{w_n} - A{y_n}} \right\|}},{\lambda _n}+p_n\} ,A{w_n} \ne A{y_n}\\
	{\lambda _n+p_n},   \quad\quad\quad\quad\quad\quad\quad\quad {\rm otherwise}.
\end{array} \right.
\end{equation*}
If ${w_n} = {y_n}$, stop and $y_n$ is a solution of the problem \eqref{1.1}. Otherwise
\item[Step 2.]Compute
 \begin{equation*}
 {x_{n + 1}} = (1 - {\theta _n}){z_n} + {\theta _n}({y_n} - {\lambda _n}(A{y_n} - A{w_n}))
\end{equation*}
Let $n= n + 1$ and return to Step 1.

\end{itemize}
\end{algorithm}

\begin{remark}
	\label{r3.2}
	\begin{itemize}
		\item[(i)] In our Algorithm \ref{alg3.1}, we can take the inertial factor $\alpha_n=1$. This is not allowed in the corresponding algorithms of \cite{CV,CH},  where the only single inertial extrapolation step is considered and the inertia
is bounded away from $1$.
\item[(ii)]In order to get the larger step sizes, being similar to the sequence $\{\theta_n\}$ in Algorithm 3.1 of \cite{W}, the sequence $\{\mu_n\}$ is used to relax the parameter $\mu$. The sequence $\{\theta_n\} $ can be called a relaxed parameter sequence, which can often improve numerical efficiency of algorithms, see \cite{CV}. If $\mu_n=0$, then the step size $\lambda_n$ is the same as the one of Algorithm 4.1 in \cite{CH}. If $\mu_n=0$ and $p_n=0$, then the step size $\lambda_n$ is the same as the one of Algorithm 3.1 of \cite{VA}.
\item[(iii)]  Note that if $\beta_n=0$, then the condition ($C_3$)(iii) can be relaxed as $0<\theta<\theta_n\le\theta_{n+1}\le \dfrac{1}{1+\varepsilon},\varepsilon  \in [0, + \infty )$, which indicates that ${\theta _n} = \hat{\theta}$ can equal 1. Setting ${\alpha _n} = \alpha$, ${\theta _n} = \hat{\theta}$, ${\mu _n} = 0$ and ${\beta _n}={p_n} = 0$, our algorithm can reduce to Algorithm 2 of \cite{CV}. In addition, if ${\alpha _n} ={\theta _n} = \alpha$, ${\mu _n} = 0$ and ${\beta _n}={p_n} = 0$,  Algorithm \ref{alg3.1} can reduce to Algorithm 1 of \cite{AB}.
		
	\end{itemize}

\end{remark}
	
\begin{lemma}
\label{l3.1} The sequences $\{\lambda_n\}$ from Algorithm 3.1 is bounded and $\lambda_n  \in [\min \{ \frac{\mu }{L},{\lambda _1}\} ,{\lambda _1} + P]$. Furthermore there exists $\lambda  \in [\min \{ \frac{\mu }{L},{\lambda _1}\} ,{\lambda _1} + P]$ such that $\mathop {\lim }\limits_{n \to \infty } {\lambda _n} = \lambda $, where $P = \sum\limits_{n = 1}^\infty  {{p_n}} $.
\end{lemma}
\begin{proof}
	By the definition of $\lambda_n$, if $A{w_n} \ne A{y_n}$, we get
	 	\begin{equation}
	 	\label{31.1}
	 	\begin{array}{l}
	 		{\lambda _n} \ge \frac{{(\mu  + {\mu _n})\left\| {w_n - y_n} \right\|}}{{\left\| {Aw_n - Ay_n} \right\|}} \ge \frac{{\mu  + {\mu _n}}}{L} \ge \frac{\mu }{L}.
	 	\end{array}
	 \end{equation}
Since  $P = \sum\limits_{n = 1}^\infty  {{p_n}} $, we have
\begin{equation}
	\begin{array}{l}
		\label{31.2}
		{\lambda _{n + 1}} \le {\lambda _n} + {p_n} \le {\lambda _1} + \sum\limits_{n = 1}^\infty  {{p_n}}={\lambda _1}+P.
	\end{array}
\end{equation}
It implies that $\{ \frac{\mu }{L},{\lambda _1}\}  \le {\lambda _n} \le {\lambda _1} + P$.\\
We have
\begin{equation}
	\label{31.4}
	\begin{array}{l}
		{\lambda _{n + 1}} - {\lambda _n} = {[{\lambda _{n + 1}} - {\lambda _n}]_ + } - {[{\lambda _{n + 1}} - {\lambda _n}]_ - }.
	\end{array}
\end{equation}
Thus
\begin{equation}
	\label{31.5}
	\begin{array}{l}
	{\lambda _{n + 1}} - {\lambda _1} = \sum\limits_{i = 1}^n {{{[{\lambda _{i + 1}} - {\lambda _i}]}_ + }}  - \sum\limits_{i = 1}^n {{{[{\lambda _{i + 1}} - {\lambda _i}]}_ - }} .
	\end{array}
\end{equation}
Since $\{\lambda_n\}$ is bounded and $\sum\limits_{n = 1}^\infty  {{{[{\lambda _{n + 1}} - {\lambda _n}]}_ + }}  \le \sum\limits_{n = 1}^\infty  {{p_n}}  < \infty $, we get $\sum\limits_{n = 1}^\infty  {{{[{\lambda _{n + 1}} - {\lambda _n}]}_ - }}$ is convergent. Therefore, there exists $\lambda  \in [\min \{ \frac{\mu }{L},{\lambda _1}\} ,{\lambda _1} + P]$ such that $\mathop {\lim }\limits_{n \to \infty } {\lambda _n} = \lambda $. This proof is completed.

\end{proof}

\begin{remark}
	\label{r3.1}  If ${w_n} = {y_n}$, then ${y_n} = {(I + {\lambda _n}B)}^{-1}(I - {\lambda _n}A){y_n}$. By Lemma \ref{l2.6}, we know ${y_n} \in \Omega$.
\end{remark}

\begin{lemma}
\label{l3.2} Suppose that the sequence $\{y_n\}$ is generated by Algorithm 3.1. Thus the following assertions hold:
\begin{itemize}
	\item[(i)]if the conditions ($C_1$) and ($C_2$) hold, then
	\begin{equation*}
		\begin{array}{l}
			{\left\| {{y_n} - {\lambda _n}(A{y_n} - A{w_n}) - p} \right\|^2} \le {\left\| {{w_n} - p} \right\|^2} - (1 - \frac{{{(\mu  + {\mu _n})^2}{\lambda _n}^2}}{{{\lambda _{n + 1}}^2}}){\left\| {{w_n} - {y_n}} \right\|^2},\forall p \in \Omega;
		\end{array}
	\end{equation*}
    \item[(ii)]if the conditions ($C_1$) and ($C_2$) hold and $A$ or $B$ is $r$-strongly monotone, then
    \begin{equation*}
    	\begin{array}{l}
    	\begin{split}
    		{\left\| {{y_n} - {\lambda _n}(A{y_n} - A{w_n}) - p} \right\|^2} &\le {\left\| {{w_n} - p} \right\|^2} - (1 - \frac{{{(\mu  + {\mu _n})^2}{\lambda _n}^2}}{{{\lambda _{n + 1}}^2}}){\left\| {{w_n} - {y_n}} \right\|^2}\\
    		&- 2r{\lambda _n}{\left\| {{y_n} - p} \right\|^2},\forall p \in \Omega.
    	\end{split}
    	\end{array}
    \end{equation*}
\end{itemize}
\end{lemma}

\begin{proof}
(i) According to the definition of $\lambda_n$,  we have that
\begin{equation}
	\label{3.1}
	\begin{array}{c}
		\begin{split}
			\left\| {{y_n} - {\lambda _n}(A{y_n} - A{w_n}) - p} \right\|^2 &= {\left\| {{y_n} - p} \right\|^2} + {\lambda _n}^2{\left\| {A{y_n} - A{w_n}} \right\|^2} - 2{\lambda _n}\left\langle {A{y_n} - A{w_n},{y_n} - p} \right\rangle \\
				&\le {\left\| {{y_n} - {w_n}} \right\|^2} + {\left\| {{w_n} - p} \right\|^2} + 2\left\langle {{y_n} - {w_n},{w_n} - p} \right\rangle \\
				&\quad+ \frac{{{(\mu  + {\mu _n})^2}{\lambda _n}^2}}{{{\lambda _{n + 1}}^2}}{\left\| {{y_n} - {w_n}} \right\|^2} - 2{\lambda _n}\left\langle {A{y_n} - A{w_n},{y_n} - p} \right\rangle \\
				&= (1 + \frac{{{(\mu  + {\mu _n})^2}{\lambda _n}^2}}{{{\lambda _{n + 1}}^2}}){\left\| {{y_n} - {w_n}} \right\|^2} + {\left\| {{w_n} - p} \right\|^2} + 2\left\langle {{y_n} - {w_n},{w_n} - {y_n}} \right\rangle \\
				&\quad+ 2\left\langle {{y_n} - {w_n},{y_n} - p} \right\rangle  - 2{\lambda _n}\left\langle {A{y_n} - A{w_n},{y_n} - p} \right\rangle \\
				&= ( - 1 + \frac{{{(\mu  + {\mu _n})^2}{\lambda _n}^2}}{{{\lambda _{n + 1}}^2}}){\left\| {{y_n} - {w_n}} \right\|^2} + {\left\| {{w_n} - p} \right\|^2}\\
				&\quad- 2\left\langle {{w_n} - {y_n} - {\lambda _n}(A{w_n} - A{y_n}),{y_n} - p} \right\rangle .
		\end{split}
	\end{array}
\end{equation}
Since ${y_n} = {(I + {\lambda _n}B)^{ - 1}}(I -{\lambda _n}A){w_n}$ and $B$ is maximal monotone, we obtain
\begin{equation}
	\label{3.2}
	\frac{{{w_n} - {y_n} - {\lambda _n}A{w_n}}}{{{\lambda _n}}} \in B{y_n}.
\end{equation}
Thus
\begin{equation*}
	\frac{{{w_n} - {y_n} - {\lambda _n}A{w_n}}}{{{\lambda _n}}} + A{y_n} \in (A + B){y_n}.
\end{equation*}
Since $A:H \to H$ is monotone and $B:H \to {2^H}$ is a maximal monotone operator, Lemma \ref{l2.3} implies that $A+B$ is maximal monotone. Since $p\in\Omega$, $0 \in (A + B)p$ and so
\begin{equation*}
	\left\langle {\frac{{{w_n} - {y_n} - {\lambda _n}A{w_n}}}{{{\lambda _n}}} + A{y_n}-0,{y_n} - p} \right\rangle  \geq 0.
\end{equation*}
Hence
\begin{equation}
\label{l+1}
	\left\langle {{w_n} - {y_n} - {\lambda _n}(A{w_n} - A{y_n}),{y_n} - p} \right\rangle \geq 0.
\end{equation}
Combing \eqref{3.1} with \eqref{l+1}, we get
\begin{equation}
	\label{3.3}
	\begin{array}{l}
		{\left\| {{y_n} - {\lambda _n}(A{y_n} - A{w_n}) - p} \right\|^2} \le {\left\| {{w_n} - p} \right\|^2} - (1 - \frac{{{(\mu  + {\mu _n})^2}{\lambda _n}^2}}{{{\lambda _{n + 1}}^2}}){\left\| {{w_n} - {y_n}} \right\|^2}.
	\end{array}
\end{equation}\\
The proof of the part (i) is completed.
 \par (ii) {\bf Case I} $B$ is $r$-strongly monotone.
 \par Since $p\in\omega$, we have $0 \in (A + B)p$ and thus $-Ap \in Bp$. Since $B$ is $r$-strongly monotone, by \eqref{3.2}, we have
 \begin{equation}
 \label{l+2}
\left\langle {\frac{{{w_n} - {y_n} - {\lambda _n}A{w_n}}}{{{\lambda _n}}} + Ap,{y_n} - p} \right\rangle  \ge r{\left\| {{y_n} - p} \right\|^2}.
\end{equation}
 The monotonicity of $A$ implies that
 \begin{equation}
 \label{l+3}
\left\langle A y_n -Ap, y_n - p \right\rangle  \geq 0.
\end{equation}
Adding together \eqref{l+2} and  \eqref{l+3}, we have
\begin{equation*}
	\left\langle {\frac{{{w_n} - {y_n} - {\lambda _n}A{w_n}}}{{{\lambda _n}}} + A{y_n},{y_n} - p} \right\rangle  \ge r{\left\| {{y_n} - p} \right\|^2},
\end{equation*}
which implies
\begin{equation}
\label{l+4}
	\left\langle {{w_n} - {y_n} - {\lambda _n}A{w_n} + {\lambda _n}A{y_n},{y_n} - p} \right\rangle  \ge r{\lambda _n}{\left\| {{y_n} - p} \right\|^2}.
\end{equation}
Utilizing \eqref{3.1} and \eqref{l+4}, we get
\begin{equation}
	\label{3.4}
	\begin{array}{l}
		\begin{split}
			{\left\| {{y_n} - {\lambda _n}(A{y_n} - A{w_n}) - p} \right\|^2} &\le {\left\| {{w_n} - p} \right\|^2} - (1 - \frac{{{(\mu  + {\mu _n})^2}{\lambda _n}^2}}{{{\lambda _{n + 1}}^2}}){\left\| {{w_n} - {y_n}} \right\|^2}\\
			&- 2r{\lambda _n}{\left\| {{y_n} - p} \right\|^2}.
		\end{split}
	\end{array}
\end{equation}
{\bf Case II} $A$ is $r$-strongly monotone.
\par The strong monotonicity of $A$ implies that
 \begin{equation}
 \label{l+5}
\left\langle A y_n -Ap, y_n - p \right\rangle  \geq r{\left\| {{y_n} - p} \right\|^2}.
\end{equation}
Since $p\in\omega$, we have $0 \in (A + B)p$ and thus $-Ap \in Bp$. Since $B$ is monotone, by \eqref{3.2}, we have
 \begin{equation}
 \label{l+6}
\left\langle {\frac{{{w_n} - {y_n} - {\lambda _n}A{w_n}}}{{{\lambda _n}}} + Ap,{y_n} - p} \right\rangle  \ge 0.
\end{equation}
The rest of proof is the same as in Case I. This completes the proof of Lemma \ref{l3.2}.
\end{proof}

\begin{lemma}
	\label{l3.3}
Assume that the conditions ($C_1$) and ($C_2$) hold,  and  $\{ {w_n}\}$ and $\{ {y_n}\} $ are sequences generated by Algorithm \ref{alg3.1}. If $\mathop {\lim }\limits_{n \to \infty } {\left\| {{w_n} - {y_n}} \right\|} = 0$ and ${\{ {w_n}\} }$ converges weakly to some $z \in H$, then $z \in \Omega $.
\end{lemma}
\begin{proof}	 Letting $(u,v) \in$Graph$(A + B)$, we get $v - Au \in Bu$. Since $B$ is maximal monotone, by \eqref{3.2}, we have
\begin{equation}
	\label{3.5}
	\left\langle v - Au - \frac{{{w_n} - {y_n} - {\lambda _n}A{w_n}}}{{{\lambda _n}}}, u - {y_n} \right\rangle  \ge 0.
\end{equation}
This implies that
\begin{equation}
	\label{3.6}
	\begin{array}{c}
			\begin{split}
		\left\langle v, u - {y_n} \right\rangle  &\ge \frac{1}{{{\lambda _n}}}\left\langle {\lambda _n}Au + {w_n} - {y_n} - {\lambda _n}A{w_n},u - {y_n} \right\rangle \\
		&= \frac{1}{{{\lambda _n}}}\left\langle {u - {y_n},{w_n} - {y_n}} \right\rangle  + \left\langle {u - {y_n},Au - A{w_n}} \right\rangle \\
		&= \frac{1}{{{\lambda _n}}}\left\langle {u - {y_n},{w_n} - {y_n}} \right\rangle  + \left\langle {u - {y_n},Au - A{y_n}} \right\rangle  + \left\langle {u - {y_n},A{y_n} - A{w_n}} \right\rangle .
			\end{split}
	\end{array}
\end{equation}
From the Lipschitz continuity and monotonicity of $A$, it follows that
\begin{equation}
	\label{3.7}
\begin{array}{c}
		\begin{split}
	\left\langle {u - {y_n},v} \right\rangle  &\ge \frac{1}{{{\lambda _n}}}\left\langle {u - {y_n},{w_n} - {y_n}} \right\rangle  + \left\langle {u - {y_n},A{y_n} - A{w_n}} \right\rangle \\
	&\ge \frac{1}{{{\lambda _n}}}\left\langle {u - {y_n},{w_n} - {y_n}} \right\rangle  - \left\| {u - {y_n}} \right\|\left\| {A{y_n} - A{w_n}} \right\|\\
	&\ge \frac{1}{{{\lambda _n}}}\left\langle {u - {y_n},{w_n} - {y_n}} \right\rangle  - L\left\| {u - {y_n}} \right\|\left\| {{y_n} - {w_n}} \right\|.
		\end{split}
\end{array}
\end{equation}
Since $\mathop {\lim }\limits_{n \to \infty } {\left\| {{w_n} - {y_n}} \right\|} = 0$, by Lemma \ref{l3.1} and \eqref{3.7}, we have
\begin{equation*}
	\left\langle {u - z,v} \right\rangle  = \mathop {\lim }\limits_{n \to \infty } \left\langle {u - {w_n},v} \right\rangle =\mathop {\lim }\limits_{n \to \infty } \left\langle {u - {y_n},v} \right\rangle  \ge 0
\end{equation*}
Since $v \in (A + B)u$ and $A+B$ is maximal monotone, we know $0 \in (A + B)z$, that is, $z \in \Omega $. This proof is completed.
\end{proof}

\begin{theorem}
\label{t3.1}Assume that the conditions ($C_1$)-($C_3$) hold. Then the sequence $\{ {x_n}\} $ from Algorithm \ref{alg3.1} converges weakly to some element $p \in \Omega $.
\end{theorem}
\begin{proof}
By the definition of $x_{n+1}$ and Lemma \ref{l2.7}, we have
\begin{equation}
	\label{3.8}
	\begin{array}{c}
		\begin{split}	
			{\left\| {{x_{n + 1}} - p} \right\|^2} &= {\left\| {(1 - {\theta _n}){z_n} + {\theta _n}({y_n} - {\lambda _n}(A{y_n} - A{w_n})) - p} \right\|^2}\\
			&= (1 - {\theta _n}){\left\| {{z_n} - p} \right\|^2} + {\theta _n}{\left\| {{y_n} - {\lambda _n}(A{y_n} - A{w_n}) - p} \right\|^2} \\&\quad- (1 - {\theta _n}){\theta _n}{\left\| {{z_n} - {y_n} + {\lambda _n}(A{y_n} - A{w_n})} \right\|^2}.
		\end{split}
	\end{array}
\end{equation}
 The definition of $x_{n+1}$ implies that  ${y_n} - {\lambda _n}(A{y_n} - A{w_n}) = \frac{{{x_{n + 1}} - (1 - {\theta _n}){z_n}}}{{{\theta _n}}}$. Thus
\begin{equation}
	\label{3.9}
	\begin{array}{c}
		\begin{split}	
			{\left\| {{z_n} - {y_n} + {\lambda _n}(A{y_n} - A{w_n})} \right\|^2} &= {\left\| {{z_n} - \frac{{{x_{n + 1}} - (1 - {\theta _n}){z_n}}}{{{\theta _n}}}} \right\|^2}\\
			&= \frac{1}{{{\theta _n^2}}}{\left\| {{x_{n + 1}} - {z_n}} \right\|^2}.
		\end{split}	
	\end{array}
\end{equation}
By Lemma \ref{l3.2}, we have
\begin{equation}
	\label{3.10}
	{\left\| {{y_n} - {\lambda _n}(A{y_n} - A{w_n}) - p} \right\|^2} \le {\left\| {{w_n} - p} \right\|^2} - (1 - \frac{{{(\mu  + {\mu _n})^2}{\lambda _n}^2}}{{{\lambda _{n + 1}}^2}}){\left\| {{w_n} - {y_n}} \right\|^2},\forall p \in \Omega .
\end{equation}
Substituting  \eqref{3.9} and \eqref{3.10} into \eqref{3.8}, we get
\begin{equation}
	\label{3.11}
	\begin{split}
		{\left\| {{x_{n + 1}} - p} \right\|^2} &\le (1 - {\theta _n}){\left\| {{z_n} - p} \right\|^2} + {\theta _n}{\left\| {{w_n} - p} \right\|^2} - {\theta _n}(1 - \frac{{{(\mu  + {\mu _n})^2}{\lambda _n}^2}}{{{\lambda _{n + 1}}^2}})\left\| {{w_n} - {y_n}} \right\|^2 \\
		&\quad- \frac{{1 - {\theta _n}}}{{{\theta _n}}}{\left\| {{x_{n + 1}} - {z_n}} \right\|^2}.
	\end{split}
\end{equation}
Lemma \ref{l3.1} and the fact$\mathop {\lim }\limits_{n \to \infty } {\mu _n} = 0$ imply that $\mathop {\lim }\limits_{n \to \infty } 1 - \frac{{{(\mu  + {\mu _n})^2}{\lambda _n}^2}}{{{\lambda _{n + 1}}^2}} = 1 - {\mu ^2} > 0$. Thus there exists a positive integer $N \ge 1$ such that
$${\theta _n}(1 - \frac{{{(\mu  + {\mu _n})^2}{\lambda _n}^2}}{{{\lambda _{n + 1}}^2}})\left\| {{w_n} - {y_n}} \right\|^2 \geq 0,~\forall~n \ge N .$$
Thanks to \eqref{3.11}, we have
\begin{equation}
	\label{3.12}
	\left\| {{x_{n + 1}} - p} \right\|^2 \le (1 - {\theta _n}){\left\| {{z_n} - p} \right\|^2} + {\theta _n}{\left\| {{w_n} - p} \right\|^2} - \frac{{1 - {\theta _n}}}{{{\theta _n}}}{\left\| {{x_{n + 1}} - {z_n}} \right\|^2},~\forall~n \ge N.
\end{equation}
By the definition of $z_n$ and Lemma \ref{l2.7}, we obtain
\begin{equation}
	\label{3.13}
	\begin{array}{c}
		\begin{split}
			{\left\| {{z_n} - p} \right\|^2} &= {\left\| {{x_n} + {\beta _n}({x_n} - {x_{n - 1}}) - p} \right\|^2}\\
			&= {\left\| {(1 + {\beta _n})({x_n} - p) - {\beta _n}({x_{n - 1}} - p)} \right\|^2}\\
			&= (1 + {\beta _n}){\left\| {{x_n} - p} \right\|^2} - {\beta _n}{\left\| {{x_{n - 1}} - p} \right\|^2} + {\beta _n}(1 + {\beta _n}){\left\| {{x_n} - {x_{n - 1}}} \right\|^2}.
		\end{split}	
	\end{array}
\end{equation}
From the definition of $w_n$ and Lemma \ref{l2.7}, it follows that
\begin{equation}
	\label{3.14}
	{\left\| {{w_n} - p} \right\|^2} = (1 + {\alpha _n}){\left\| {{x_n} - p} \right\|^2} - {\alpha _n}{\left\| {{x_{n - 1}} - p} \right\|^2} + {\alpha _n}(1 + {\alpha _n}){\left\| {{x_n} - {x_{n - 1}}} \right\|^2}.
\end{equation}
The definition of $z_n$ means that
\begin{equation}
	\label{3.15}
	\begin{array}{c}
		\begin{split}
			{\left\| {{x_{n + 1}} - {z_n}} \right\|^2} &= {\left\| {{x_{n + 1}} - {x_n} - {\beta _n}({x_n} - {x_{n - 1}})} \right\|^2}\\
			&\ge {\left\| {{x_{n + 1}} - {x_n}} \right\|^2} + {\beta _n}^2{\left\| {{x_n} - {x_{n - 1}}} \right\|^2} - 2{\beta _n}\left\| {{x_{n + 1}} - {x_n}} \right\|\left\| {{x_n} - {x_{n - 1}}} \right\|\\
			&\ge (1 - {\beta _n}){\left\| {{x_{n + 1}} - {x_n}} \right\|^2} + ({\beta _n}^2 - {\beta _n}){\left\| {{x_n} - {x_{n - 1}}} \right\|^2}.
		\end{split}
	\end{array}
\end{equation}
Owing to \eqref{3.13}, \eqref{3.14}, \eqref{3.15} and \eqref{3.12}, we know
\begin{equation}
	\label{3.16}
	\begin{split}
		\left\| {{x_{n + 1}} - p} \right\|^2 &\le (1 - {\theta _n})[(1 + {\beta _n}){\left\| {{x_n} - p} \right\|^2} - {\beta _n}{\left\| {{x_{n - 1}} - p} \right\|^2} + {\beta _n}(1 + {\beta _n}){\left\| {{x_n} - {x_{n - 1}}} \right\|^2}]\\
		&\quad 	 + {\theta _n}[(1 + {\alpha _n}){\left\| {{x_n} - p} \right\|^2} - {\alpha _n}{\left\| {{x_{n - 1}} - p} \right\|^2} + {\alpha _n}(1 + {\alpha _n}){\left\| {{x_n} - {x_{n - 1}}} \right\|^2}]\\
		&\quad - \frac{{1 - {\theta _n}}}{{{\theta _n}}}[(1 - {\beta _n}){\left\| {{x_{n + 1}} - {x_n}} \right\|^2} + ({\beta _n}^2 - {\beta _n}){\left\| {{x_n} - {x_{n - 1}}} \right\|^2}] \\
		&\le [1 + (1 - {\theta _n}){\beta _n} + {\theta _n}{\alpha _n}]{\left\| {{x_n} - p} \right\|^2} - [(1 - {\theta _n}){\beta _n} + {\theta _n}{\alpha _n}]{\left\| {{x_{n - 1}} - p} \right\|^2}\\
		&\quad+ [(1 - {\theta _n})(1 + {\beta _n}){\beta _n} + {\theta _n}(1 + {\alpha _n}){\alpha _n} - \frac{{(1 - {\theta _n})}}{{{\theta _n}}}({\beta _n}^2 - {\beta _n})]{\left\| {{x_n} - {x_{n - 1}}} \right\|^2}\\
		&\quad- \frac{{(1 - {\theta _n})}}{{{\theta _n}}}(1 - {\beta _n}){\left\| {{x_{n + 1}} - {x_n}} \right\|^2} \\
		& \le (1 + {a_n}){\left\| {{x_n} - p} \right\|^2} - {a_n}{\left\| {{x_{n - 1}} - p} \right\|^2} + {b_n}{\left\| {{x_n} - {x_{n - 1}}} \right\|^2} - {c_n}{\left\| {{x_{n + 1}} - {x_n}} \right\|^2}.
	\end{split}
\end{equation}
where ${a_n} = (1 - {\theta _n}){\beta _n} + {\theta _n}{\alpha _n}$, ${b_n} := (1 - {\theta _n})(1 + {\beta _n}){\beta _n} + {\theta _n}(1 + {\alpha _n}){\alpha _n} - \frac{{1 - {\theta _n}}}{{{\theta _n}}}({\beta _n}^2 - {\beta _n})$ and ${c_n} := \frac{{1 - {\theta _n}}}{{{\theta _n}}}(1 - {\beta _n})$. \\
Define
\begin{equation*}
	{\Gamma _n}: = {\left\| {{x_n} - p} \right\|^2} - {a_n}{\left\| {{x_{n - 1}} - p} \right\|^2} + {b_n}{\left\| {{x_n} - {x_{n - 1}}} \right\|^2}
\end{equation*}
Since $\{a_n\}$ is non-decreasing, by \eqref{3.16}, we deduce that
\begin{equation}
	\label{3.17}
	\begin{array}{c}
		\begin{split}
			{\Gamma _{n + 1}} - {\Gamma _n} &= {\left\| {{x_{n + 1}} - p} \right\|^2} - (1 + {a_{n + 1}}){\left\| {{x_n} - p} \right\|^2} + {a_n}{\left\| {{x_{n - 1}} - p} \right\|^2} \\
			&\quad- {b_n}{\left\| {{x_n} - {x_{n - 1}}} \right\|^2} + {b_{n + 1}}{\left\| {{x_{n + 1}} - {x_n}} \right\|^2}\\
			&\le {\left\| {{x_{n + 1}} - p} \right\|^2} - (1 + {a_n}){\left\| {{x_n} - p} \right\|^2} + {a_n}{\left\| {{x_{n - 1}} - p} \right\|^2} \\
			&\quad- {b_n}{\left\| {{x_n} - {x_{n - 1}}} \right\|^2} + {b_{n + 1}}{\left\| {{x_{n + 1}} - {x_n}} \right\|^2}\\
			&\le  - ({c_n} - {b_{n + 1}}){\left\| {{x_{n + 1}} - {x_n}} \right\|^2}.
		\end{split}
	\end{array}
\end{equation}
The condition ($C_3$) means that
\begin{equation}
	\label{3.18}
	\begin{array}{c}
		\begin{split}
			{c_n} - {b_{n + 1}} &= \frac{{1 - {\theta _n}}}{{{\theta _n}}}(1 - {\beta _n}) - (1 - {\theta _{n + 1}})(1 + {\beta _{n + 1}}){\beta _{n + 1}} - {\theta _{n + 1}}(1 + {\alpha _{n + 1}}){\alpha _{n + 1}}\\
			&\quad+ \frac{{1 - {\theta _{n + 1}}}}{{{\theta _{n + 1}}}}({\beta _{n + 1}}^2 - {\beta _{n + 1}})\\
			&\ge \frac{{1 - {\theta _{n + 1}}}}{{{\theta _{n + 1}}}}(1 - 2{\beta _{n + 1}} + {\beta _{n + 1}}^2)  - (1 - {\theta _{n + 1}})({\beta _{n + 1}} + {\beta _{n + 1}}^2)- 2{\theta _{n + 1}}\\
			&\ge \varepsilon (1 - 2\beta  + {\beta ^2})  - \frac{\varepsilon }{{1 + \varepsilon }}(\beta  + {\beta ^2})- \frac{2}{{1 + \varepsilon }}\\
			&= \frac{1}{{1 + \varepsilon }}[{\varepsilon ^2}{\beta ^2} - (3\varepsilon  + 2{\varepsilon ^2})\beta  + ({\varepsilon ^2} + \varepsilon  - 2)].
		\end{split}
	\end{array}
\end{equation}
Combining \eqref{3.17} and \eqref{3.18}, we infer that
\begin{equation}
	\label{3.19}
	{\Gamma _{n + 1}} - {\Gamma _n} \le  - \delta {\left\| {{x_{n + 1}} - {x_n}} \right\|^2}.
\end{equation}
where $\delta : = \frac{{{\varepsilon ^2}{\beta ^2} - (3\varepsilon  + 2{\varepsilon ^2})\beta  + ({\varepsilon ^2} + \varepsilon  - 2)}}{{1 + \varepsilon }}$. Since $\beta<\frac{{3 + 2\varepsilon  - \sqrt {8\varepsilon  + 17} }}{{2\varepsilon }},\varepsilon  \in (1, + \infty )$ , we conclude that $ \delta >0$.  From \eqref{3.19}, it follows that
\begin{equation*}
	{\Gamma _{n + 1}} - {\Gamma _n} \le 0.
\end{equation*}
Thus the sequence $\{ {\Gamma _n}\}$ is nonincreasing. The condition ($C_3$) implies that $b_n>0$. Thus
\begin{equation}
	\label{3.20}
	\begin{array}{c}
		\begin{split}
			{\Gamma _n} &= {\left\| {{x_n} - p} \right\|^2} - {a_n}{\left\| {{x_{n - 1}} - p} \right\|^2} + {b_n}{\left\| {{x_n} - {x_{n - 1}}} \right\|^2}\\
			& \ge {\left\| {{x_n} - p} \right\|^2} - {a_n}{\left\| {{x_{n - 1}} - p} \right\|^2}.
		\end{split}
	\end{array}
\end{equation}
This implies that
\begin{equation}
	\label{3.21}
	\begin{array}{c}
		\begin{split}
			{\left\| {{x_n} - p} \right\|^2} &\le {a_n}{\left\| {{x_{n - 1}} - p} \right\|^2} + {\Gamma _n}\\
			&\le a{\left\| {{x_{n - 1}} - p} \right\|^2} + {\Gamma _1}\\
			&\vdots \\
			&\le {a^n}{\left\| {{x_0} - p} \right\|^2} + {\Gamma _1}(1 + a +  \ldots {a^{n - 1}})\\
			&\le {a^n}{\left\| {{x_0} - p} \right\|^2} + \frac{{{\Gamma _1}}}{{1 - a}}.
		\end{split}
	\end{array}
\end{equation}
where $a :=\frac{{5 + 2\varepsilon  - \sqrt {8\varepsilon  + 17} }}{{2 + 2\varepsilon }} < 1$. \\By definition of $\{ {\Gamma _n}\}$, we have
\begin{equation}
	\label{3.22}
	\begin{array}{c}
		\begin{split}
			{\Gamma _{n + 1}} &= {\left\| {{x_{n + 1}} - p} \right\|^2} - {a_{n + 1}}{\left\| {{x_n} - p} \right\|^2} + {b_{n + 1}}{\left\| {{x_{n + 1}} - {x_n}} \right\|^2}\\
			&\ge  - {a_{n + 1}}{\left\| {{x_n} - p} \right\|^2}.
		\end{split}
	\end{array}
\end{equation}
According to \eqref{3.21} and \eqref{3.22}, we infer that
\begin{equation}
	\label{3.23}
	- {\Gamma _{n + 1}} \le {a_{n + 1}}{\left\| {{x_n} - p} \right\|^2} \le a{\left\| {{x_n} - p} \right\|^2} \le {a^{n + 1}}{\left\| {{x_0} - p} \right\|^2} + \frac{{a{\Gamma _1}}}{{1 - a}}.
\end{equation}
By \eqref{3.23} and \eqref{3.19}, we get
\begin{equation}
	\label{3.24}
	\begin{array}{c}
		\begin{split}
			\delta \sum\limits_{k = 1}^n {{{\left\| {{x_{k + 1}} - {x_k}} \right\|}^2} \le {\Gamma _1}}  - {\Gamma _{n + 1}} &\le {a^{n + 1}}{\left\| {{x_0} - p} \right\|^2} + \frac{{{\Gamma _1}}}{{1 - a}}\\
			&\le {\left\| {{x_0} - p} \right\|^2} + \frac{{{\Gamma _1}}}{{1 - a}}.
		\end{split}
	\end{array}
\end{equation}
This implies \begin{equation}
	\label{l+1.1}
\sum \limits_{k = 1}^\infty \| x_{k + 1}- x_k \|^2 \leq \infty .
\end{equation}
 Thus
\begin{equation}
	\label{3.25}
	\left\| {{x_{n + 1}} - {x_n}} \right\| \to 0,n \to \infty.
\end{equation}
The definition of $\{ {\Gamma _n}\}$ implies that
\begin{equation}
	\label{3.26}
	{\left\| {{x_{n + 1}} - {w_n}} \right\|^2} = {\left\| {{x_{n + 1}} - {x_n}} \right\|^2} + {\alpha _n}{\left\| {{x_n} - {x_{n - 1}}} \right\|^2} - 2{\alpha _n}\left\langle {{x_{n + 1}} - {x_n},{x_n} - {x_{n - 1}}} \right\rangle .
\end{equation}
Since $\{\alpha_n\}$ is bounded, by \eqref{3.26}, we obtain
\begin{equation}
		\label{3.27}
	{\left\| {{x_{n + 1}} - {w_n}} \right\|^2} \to 0,n \to \infty .
\end{equation}
On the other hand
\begin{equation}
	\label{3.28}
	\left\| {{x_n} - {w_n}} \right\| \le \left\| {{x_n} - {x_{n + 1}}} \right\| + \left\| {{x_{n + 1}} - {w_n}} \right\|.
\end{equation}
From \eqref{3.25} and \eqref{3.27}, it follows that
\begin{equation}
		\label{3.29}
	{\left\| {{x_n} - {w_n}} \right\|^2} \to 0,n \to \infty .
\end{equation}
By \eqref{3.16}, we have
\begin{equation}
	\label{3.30}
	\begin{array}{c}
		\begin{split}
			{\left\| {{x_{n+1}} - p} \right\|^2} &\le (1 + {a_n}){\left\| {{x_n} - p} \right\|^2} - {a_n}{\left\| {{x_{n - 1}} - p} \right\|^2} + {b_n}{\left\| {{x_n} - {x_{n - 1}}} \right\|^2} - {c_n}{\left\| {{x_{n + 1}} - {x_n}} \right\|^2}\\
			&\le (1 + {a_n}){\left\| {{x_n} - p} \right\|^2} - {a_n}{\left\| {{x_{n - 1}} - p} \right\|^2} + {b_n}{\left\| {{x_n} - {x_{n - 1}}} \right\|^2}\\
			&= {\left\| {{x_n} - p} \right\|^2} + {a_n}({\left\| {{x_n} - p} \right\|^2}-{\left\| {{x_{n - 1}} - p} \right\|^2}) + {b_n}{\left\| {{x_n} - {x_{n - 1}}} \right\|^2}.
		\end{split}
	\end{array}
\end{equation}
Since $0 \le a_n < a <1$ and $\{b_n\}$ is bounded, by Lemma \ref{l2.1} and \eqref{l+1.1}, we know there exists $l\in [0, +\infty)$ such that
\begin{equation}
	\label{3.31}
	\mathop {\lim }\limits_{n \to \infty } {\left\| {{x_n} - p} \right\|^2} = l.
\end{equation}
Then from \eqref{3.14}, we have
\begin{equation}
	\label{3.32}
		\begin{array}{c}
		\begin{split}
	{\left\| {{w_n} - p} \right\|^2} &=(1 + {\alpha _n}){\left\| {{x_n} - p} \right\|^2} - {\alpha _n}{\left\| {{x_{n - 1}} - p} \right\|^2} + {\alpha _n}(1 + {\alpha _n}){\left\| {{x_n} - {x_{n - 1}}} \right\|^2}\\
	&= {\left\| {{x_n} - p} \right\|^2} + {\alpha _n}({\left\| {{x_n} - p} \right\|^2} - {\left\| {{x_{n - 1}} - p} \right\|^2}) + {\alpha _n}(1 + {\alpha _n}){\left\| {{x_n} - {x_{n - 1}}} \right\|^2}.
			\end{split}
\end{array}
\end{equation}
Since $\{ {\alpha _n}\} $ is bounded, by \eqref{3.25} and \eqref{3.31}, we obtain
\begin{equation}
	\label{3.33}
	\mathop {\lim }\limits_{n \to \infty } {\left\| {{w_n} - p} \right\|^2} = l.
\end{equation}
 Utilizing the similar discussion as in obtaining \eqref{3.33}, we can get
\begin{equation}
	\label{3.34}
	\mathop {\lim }\limits_{n \to \infty } {\left\| {{z_n} - p} \right\|^2} = l.
\end{equation}
Owing to \eqref{3.11}, we have
\begin{equation}
	\label{3.35}
	\begin{array}{c}
		\begin{split}	
			{\theta _n}(1 - \frac{{{(\mu  + {\mu _n})^2}{\lambda _n}^2}}{{{\lambda _{n + 1}}^2}}){\left\| {{w_n} - {y_n}} \right\|^2} &\le (1 - {\theta _n}){\left\| {{z_n} - p} \right\|^2} + {\theta _n}{\left\| {{w_n} - p} \right\|^2} \\
			&\quad- \frac{{(1 - {\theta _n})}}{{{\theta _n}}}{\left\| {{x_{n + 1}} - {z_n}} \right\|^2} - {\left\| {{x_{n + 1}} - p} \right\|^2}\\
			&\le (1 - {\theta _n}){\left\| {{z_n} - p} \right\|^2} + {\theta _n}{\left\| {{w_n} - p} \right\|^2} - {\left\| {{x_{n + 1}} - p} \right\|^2}.
		\end{split}
	\end{array}
\end{equation}
Since $0<\theta_n<\dfrac{1}{1+\epsilon},\varepsilon  \in (1, + \infty )$ and $\mathop {\lim }\limits_{n \to \infty } 1 - \frac{{{(\mu  + {\mu _n})^2}{\lambda _n}^2}}{{{\lambda _{n + 1}}^2}} = 1 - {\mu ^2} > 0$, by \eqref{3.31} \eqref{3.33} and \eqref{3.34}, we get $\mathop {\lim }\limits_{n \to \infty } \left\| {{w_n} - {y_n}} \right\| = 0. $

Since $\{ {x_n}\} $ is bounded, we assume that there exists a subsequence $\{ {x_{{n_k}}}\} $ of $\{ {x_n}\} $
such that  ${x_{{n_k}}}\rightharpoonup z \in H$. The fact $\lim_{n\to\infty}\| {{w_n} - {x_n}}\|= 0$ implies that ${w_{{n_k}}}\rightharpoonup z \in H$. By Lemma \ref{l3.3}, we know $z \in \Omega $.
The two assumptions of Lemma \ref{l2.2} are verified.  Lemma \ref{l2.2} ensure that the sequence $\{ {x_n}\} $ converges weakly to  $\mu^* \in \Omega $. The proof is completed.
\end{proof}
\par Let $C$ is a nonempty, closed and convex subset of $H$. If $B=N_C$, then MIP (\ref{1.1}) reduce to the following variational inequality, denoted by VI$(A,C)$: find a point $x^*\in C$ such that
\begin{equation*}
\langle A(x^*), y-x^*\rangle\geq 0,~\forall~ y\in C.
\end{equation*}
Denote the solution set of VI$(A,C)$ by $S$.
\begin{assumption}
	\label{a3.1}
	\begin{itemize}
	\item[(i)] The solution set of the problem (VI) is nonempty, that is, $S \ne \emptyset $.
	\item[(ii)] The mapping $A:H\to H$ is pseudomonotone, Lipschitz continuous and $A$ satisfies the condition: for any $\{x_n\}\subset H$ with $x_n\rightharpoonup w^*$, one has $\|Aw^*\|\leq \liminf_{n\to \infty}\|Ax_n\|.$
	\end{itemize}
\end{assumption}

\begin{proposition}
	\label{p3.1} Suppose that $B=N_C$ in Algorithm 3.1. Thus the following statements hold:
	\begin{itemize}
		\item[(i)]if the Assumption \ref{a3.1} hold, then
		\begin{equation*}
			\begin{array}{l}
				{\left\| {{y_n} - {\lambda _n}(A{y_n} - A{w_n}) - p} \right\|^2} \le {\left\| {{w_n} - p} \right\|^2} - (1 - \frac{{{(\mu  + {\mu _n})^2}{\lambda _n}^2}}{{{\lambda _{n + 1}}^2}}){\left\| {{w_n} - {y_n}} \right\|^2},\forall p \in S;
			\end{array}
		\end{equation*}
		\item[(ii)]if the Assumption \ref{a3.1} hold and $A$ is $\mu$-strongly pseudomonotone, then
		\begin{equation*}
			\begin{array}{l}
				\begin{split}
					{\left\| {{y_n} - {\lambda _n}(A{y_n} - A{w_n}) - p} \right\|^2} &\le {\left\| {{w_n} - p} \right\|^2} - (1 - \frac{{{(\mu  + {\mu _n})^2}{\lambda _n}^2}}{{{\lambda _{n + 1}}^2}}){\left\| {{w_n} - {y_n}} \right\|^2}\\
					&- 2\mu{\lambda _n}{\left\| {{y_n} - p} \right\|^2},\forall p \in S.
				\end{split}
			\end{array}
		\end{equation*}
	\end{itemize}
\end{proposition}
\begin{proof}
	Since $B=N_C$, we know ${y_n} = {(I + {\lambda _n}B)}^{-1}(I - {\lambda _n}A){w_n} = {P_C}({w_n} - {\lambda _n}A{w_n})$ . Thus we get
	\begin{equation}
		\label{c32}
		\left\langle {{y_n} - {w_n} + {\lambda _n}A{w_n},{y_n} - p} \right\rangle  \le 0.
	\end{equation}
	For any given $p\in S$, $\left\langle {Ap,y_n - p} \right\rangle  \ge 0.$  Since $A$ is  pseudomonotone, we have
	\begin{equation}
		\label{c31}
		\left\langle {Ay_n,y_n - p} \right\rangle  \ge 0.
	\end{equation}
	From \eqref{c32} and \eqref{c31}, it follows that
	\begin{equation}
		\label{c33}
		\left\langle {{w_n} - {y_n} - {\lambda _n}A{w_n} + {\lambda _n}A{y_n},{y_n} - p} \right\rangle  \ge 0.
	\end{equation}
	which is just \eqref{l+1}. The rest of the proof follows the same arguments as in Lemma \ref{l3.2}. The proof of the part (i) is completed.
 \par (ii)
For any given $p\in S$, $\left\langle {Ap,y_n - p} \right\rangle  \ge 0.$ The strong pseudomonotonicity of $A$ implies that
 \begin{equation}
 \label{l2+5}
\left\langle A y_n, y_n - p \right\rangle  \geq \mu{\left\| {{y_n} - p} \right\|^2}.
\end{equation}
According to\eqref{c32} and \eqref{l2+5}, we have
 \begin{equation}
 \label{l2+7}
\langle w_n-y_n-\lambda_n(Aw_n-Ay_n), y_n - p \rangle \geq \mu\lambda_n{\left\| {{y_n} - p} \right\|^2}.
\end{equation}
The rest of proof is the same as in Lemma \ref{l3.2}. This completes the proof of Proposition \ref{p3.1}.
\end{proof}

\begin{proposition}
	\label{p3.2} Assume that $B=N_C$, the Assumption \ref{a3.1} hold,  and  $\{ {w_n}\}$ and $\{ {y_n}\} $ are sequences generated by Algorithm \ref{alg3.1}. If $\mathop {\lim }\limits_{n \to \infty } {\left\| {{w_n} - {y_n}} \right\|} = 0$ and ${\{ {w_n}\} }$ converges weakly to some $z \in H$, then $z \in S $.
\end{proposition}
\begin{proof}
	This proof is the same as in Lemma 3.7 of \cite{TYCR}, and we omit it.
\end{proof}

\begin{corollary}
	\label{c3.3} If $B=N_C$, the conditions ($C_3$) and Assumption \ref{a3.1} hold, then the sequence $\{ {x_n}\} $ from Algorithm \ref{alg3.1} converges weakly to a point $p\in S$.
\end{corollary}
\begin{proof}
Replacing Lemmas \ref{l3.2} and \ref{l3.3} by  Propositions  \ref{p3.1} and \ref{p3.2}, respectively and using the same proof as in Theorem \ref{t3.1}, we obtain the desired conclusion.
\end{proof}

\begin{remark}
	\label{r3.3}
  Compared with Theorem 4.2 of \cite{YI}, the advantages of Corollary \ref{c3.3} have (i) the sequence $\{\alpha_n\}$ may not be non-decreasing;  (ii) the sequence $\{\beta_n\}$ may not be a constant; (iii) we require $0<\theta<\theta_n\le\theta_{n+1}<\dfrac{1}{1+\epsilon},\varepsilon  \in (1, + \infty )$ other than $ \varepsilon  \in (2, + \infty )$, which extend the taking value interval of $\theta_n$.  From the numerical experiment in Section 5, it can be seen that the larger the values of $\theta_{n}$, the better the algorithm performs.
\end{remark}
 \par Motivated by the Theorem 5.1 of \cite{SILD}, which may be the first nonasymptotic convergence rate results of inertial projection-type algorithm for solving variational inequalities with monotone mappings, we give the nonasymptotic $O(\frac{1}{\sqrt{n}})$ convergence rate with $``\min_{N\leq i \leq n}"$ for Algorithm \ref{alg3.1}.
\begin{theorem}	
	\label{t3.2}Assume that the conditions ($C_1$)-($C_3$) hold, the sequence $\{x_n\}$ be generated by Algorithm \ref{alg3.1} and $[t]_+ := max\{t, 0\}$. Then for any $p \in \Omega $, there exist constants M  such that the following estimate holds
	\begin{equation*}
	\mathop {\min }\limits_{N \le i \le n} \left\| {{w_i} - {y_i}} \right\| \le {(\frac{{({{\left\| {{x_N} - p} \right\|}^2} + \frac{a}{{1 - a}}{{[{{\left\| {{x_N} - p} \right\|}^2} - {{\left\| {{x_{N - 1}} - p} \right\|}^2}]}_ + } + \frac{{(2 + \frac{{1 - \theta }}{{4\theta }})M}}{{1 - a}})\frac{1}{{\theta (1 - \mu )}}}}{{n - N + 1}})^{\frac{1}{2}}}.
	\end{equation*}
\end{theorem}	
\begin{proof}
	Lemma \ref{l3.1} and the fact$\mathop {\lim }\limits_{n \to \infty } {\mu _n} = 0$ imply that $\mathop {\lim }\limits_{n \to \infty } 1 - \frac{{{(\mu  + {\mu _n})^2}{\lambda _n}^2}}{{{\lambda _{n + 1}}^2}} = 1 - {\mu ^2} > 0$. Thus there exists a positive integer $N \ge 1$ such that
	$${\theta _n}(1 - \frac{{{(\mu  + {\mu _n})^2}{\lambda _n}^2}}{{{\lambda _{n + 1}}^2}})\left\| {{w_n} - {y_n}} \right\|^2 \geq \theta _n(1 -\mu)\left\| {{w_n} - {y_n}} \right\|^2,~\forall~n \ge N .$$
	Thanks to \eqref{3.11}, we have, for all $n \ge N$
	\begin{equation}
		\label{+1}
		\begin{split}
			{\left\| {{x_{n + 1}} - p} \right\|^2} &\le (1 - {\theta _n}){\left\| {{z_n} - p} \right\|^2} + {\theta _n}{\left\| {{w_n} - p} \right\|^2} - {\theta _n}(1-\mu)\left\| {{w_n} - {y_n}} \right\|^2 - \frac{{1 - {\theta _n}}}{{{\theta _n}}}{\left\| {{x_{n + 1}} - {z_n}} \right\|^2}\\
			&\le (1 - {\theta _n}){\left\| {{z_n} - p} \right\|^2} + {\theta _n}{\left\| {{w_n} - p} \right\|^2} - {\theta}(1-\mu)\left\| {{w_n} - {y_n}} \right\|^2 - \frac{{1 - {\theta _n}}}{{{\theta _n}}}{\left\| {{x_{n + 1}} - {z_n}} \right\|^2}.
		\end{split}
	\end{equation}
Owing to \eqref{3.13}, \eqref{3.14}, \eqref{3.15} and \eqref{+1}, we know, for all $n \ge N$
\begin{equation}
	\label{+2}
	\begin{split}
		\left\| {{x_{n + 1}} - p} \right\|^2 &\le (1 - {\theta _n})[(1 + {\beta _n}){\left\| {{x_n} - p} \right\|^2} - {\beta _n}{\left\| {{x_{n - 1}} - p} \right\|^2} + {\beta _n}(1 + {\beta _n}){\left\| {{x_n} - {x_{n - 1}}} \right\|^2}]\\
		&\quad 	 + {\theta _n}[(1 + {\alpha _n}){\left\| {{x_n} - p} \right\|^2} - {\alpha _n}{\left\| {{x_{n - 1}} - p} \right\|^2} + {\alpha _n}(1 + {\alpha _n}){\left\| {{x_n} - {x_{n - 1}}} \right\|^2}]\\
		&\quad - \frac{{1 - {\theta _n}}}{{{\theta _n}}}[(1 - {\beta _n}){\left\| {{x_{n + 1}} - {x_n}} \right\|^2} + ({\beta _n}^2 - {\beta _n}){\left\| {{x_n} - {x_{n - 1}}} \right\|^2}] \\
		&\quad- {\theta}(1-\mu)\left\| {{w_n} - {y_n}} \right\|^2\\
		&\le [1 + (1 - {\theta _n}){\beta _n} + {\theta _n}{\alpha _n}]{\left\| {{x_n} - p} \right\|^2} - [(1 - {\theta _n}){\beta _n} + {\theta _n}{\alpha _n}]{\left\| {{x_{n - 1}} - p} \right\|^2}\\
		&\quad+ [(1 - {\theta _n})(1 + {\beta _n}){\beta _n} + {\theta _n}(1 + {\alpha _n}){\alpha _n} - \frac{{(1 - {\theta _n})}}{{{\theta _n}}}({\beta _n}^2 - {\beta _n})]{\left\| {{x_n} - {x_{n - 1}}} \right\|^2}\\
		&\quad- \frac{{(1 - {\theta _n})}}{{{\theta _n}}}(1 - {\beta _n}){\left\| {{x_{n + 1}} - {x_n}} \right\|^2} - {\theta}(1-\mu)\left\| {{w_n} - {y_n}} \right\|^2\\
		& \le (1 + {a_n}){\left\| {{x_n} - p} \right\|^2} - {a_n}{\left\| {{x_{n - 1}} - p} \right\|^2} + {b_n}{\left\| {{x_n} - {x_{n - 1}}} \right\|^2} - {c_n}{\left\| {{x_{n + 1}} - {x_n}} \right\|^2}\\
		&\quad - {\theta}(1-\mu)\left\| {{w_n} - {y_n}} \right\|^2.\\
		& \le (1 + {a_n}){\left\| {{x_n} - p} \right\|^2} - {a_n}{\left\| {{x_{n - 1}} - p} \right\|^2} + {b_n}{\left\| {{x_n} - {x_{n - 1}}} \right\|^2}- {\theta}(1-\mu)\left\| {{w_n} - {y_n}} \right\|^2.
	\end{split}
\end{equation}
where $\{a_n\}$, $\{b_n\}$ and $\{c_n\}$ have the same definitions as in \eqref{3.16}. \\

This implies that, for all $n \ge N$
\begin{equation}
	\label{+3}
		 {\theta}(1-\mu)\left\| {{w_n} - {y_n}} \right\|^2 \le {\left\| {{x_n} - p} \right\|^2} -\left\| {{x_{n + 1}} - p} \right\|^2+{a_n}({\left\| {{x_n} - p} \right\|^2}-{\left\| {{x_{n - 1}} - p} \right\|^2}) + {b_n}{\left\| {{x_n} - {x_{n - 1}}} \right\|^2} .
\end{equation}
Let ${\sigma _n} = {\left\| {{x_n} - p} \right\|^2},{K_n} = {\sigma _n} - {\sigma _{n - 1}}$ and ${\tau _n} ={b_n}{\left\| {{x_n} - {x_{n - 1}}} \right\|^2} $. We get, for all $n \ge N$
\begin{equation}
	\label{+4}
		\begin{split}
	{\theta}(1-\mu)\left\| {{w_n} - {y_n}} \right\|^2 &\le \sigma_n-\sigma_{n+1}+a_nK_n+\tau_n\\
	&\le \sigma_n-\sigma_{n+1}+a_n[K_n]_++\tau_n\\
	&\le \sigma_n-\sigma_{n+1}+a[K_n]_++\tau_n\\
		\end{split}
\end{equation}
where $a$ has been defined in \eqref{3.21}.\\
 In view of \eqref{l+1.1}, we have $\sum \limits_{k = 1}^\infty \| x_{k + 1}- x_k \|^2 \leq \infty$. Thus, there exists a positive constant $M$ such that \[\sum \limits_{k = 1}^\infty \| x_{k + 1}- x_k \|^2 \leq M\]
Therefore,
\begin{equation}
	\label{+5}
	\begin{split}
			\sum\limits_{n = 1}^\infty  {{\tau _n}}  &= \sum\limits_{n = 1}^\infty  {{b_n}{{\left\| {{x_n} - {x_{n - 1}}} \right\|}^2}} \\
			&= \sum\limits_{n = 1}^\infty  {[(1 - {\theta _n})(1 + {\beta _n}){\beta _n} + {\theta _n}(1 + {\alpha _n}){\alpha _n} + \frac{{1 - {\theta _n}}}{{{\theta _n}}}({\beta _n} - {\beta _n}^2)]{{\left\| {{x_n} - {x_{n - 1}}} \right\|}^2}} \\
			&\le \sum\limits_{n = 1}^\infty  {(2 + \frac{{1 - \theta }}{{4\theta }}){{\left\| {{x_n} - {x_{n - 1}}} \right\|}^2}} \\
			&= (2 + \frac{{1 - \theta }}{{4\theta }})\sum\limits_{n = 1}^\infty  {{{\left\| {{x_n} - {x_{n - 1}}} \right\|}^2}} \\
			&\le (2 + \frac{{1 - \theta }}{{4\theta }})M = {C_1}
	\end{split}
\end{equation}
From \eqref{+2}, it follows that, for all $n \ge N$
\begin{equation*}
	{K_{n + 1}} \le {a_n}{K_n} + {\tau _n} \le {a_n}{[{K_n}]_ + } + {\tau _n}.
\end{equation*}
Thus,
\begin{equation}
	\label{+6}
	\begin{split}
		{[{K_{n + 1}}]_ + } &\le a{[{K_n}]_ + } + {\tau _n}\\
		&\le {a^{n - N + 1}}{[{K_n}]_ + } + \sum\limits_{j = 1}^{n - N + 1} {{a^{j - 1}}{\tau _{n + 1 - j}}}.
	\end{split}
\end{equation}
Combining  \eqref{+5} and \eqref{+6}, we obtain
\begin{equation}
	\label{+7}
	\begin{split}
		\sum\limits_{n = N}^\infty  {{{[{K_{n + 1}}]}_ + }}  &\le \frac{a}{{1 - a}}{[{K_N}]_ + } + \frac{1}{{1 - a}}\sum\limits_{n = N}^\infty  {{\tau _n}} \\
		&\le \frac{a}{{1 - a}}{[{K_N}]_ + } + \frac{{{C_1}}}{{1 - a}}.
	\end{split}
\end{equation}
\end{proof}
From \eqref{+4}, it follows that
\begin{equation}
	\label{+8}
	\begin{split}
		\theta (1 - \mu )\sum\limits_{n = N}^n {{{\left\| {{w_n} - {y_n}} \right\|}^2}}  &\le {\sigma _N} - \sigma _{{n + 1}} + a\sum\limits_{n = N}^n {{{[{K_n}]}_ + }}  + \sum\limits_{n = N}^n {{\tau _n}} \\
		&\le {\sigma _N} + a({[{K_N}]_ + } + \sum\limits_{n = N}^n {{{[{K_{n + 1}}]}_ + }} ) + \sum\limits_{n = N}^n {{\tau _n}} \\
		&\le {\sigma _N} + a{[{K_N}]_ + } + \frac{{{a^2}}}{{1 - a}}{[{K_N}]_ + } + \frac{{a{C_1}}}{{1 - a}} + {C_1}\\
		&= {\sigma _N} + \frac{a}{{1 - a}}{[{K_N}]_ + } + \frac{{{C_1}}}{{1 - a}}\\
		&={\sigma _N} + \frac{a}{{1 - a}}{[{K_N}]_ + } + \frac{{(2 + \frac{{1 - \theta }}{{4\theta }})M}}{{1 - a}}.
	\end{split}
\end{equation}
This implies that
\begin{equation}
	\label{+9}
     \sum\limits_{i = N}^n {{{\left\| {{w_i} - {y_i}} \right\|}^2}}  \le ({\left\| {{x_N} - p} \right\|^2} + \frac{a}{{1 - a}}{[{\left\| {{x_N} - p} \right\|^2} - {\left\| {{x_{N - 1}} - p} \right\|^2}]_ + } + \frac{{(2 + \frac{{1 - \theta }}{{4\theta }})M}}{{1 - a}})\frac{1}{{\theta (1 - \mu )}}.
\end{equation}
and thus
\begin{equation}
	\label{+10}
	\mathop {\min }\limits_{N \le i \le n} {\left\| {{w_i} - {y_i}} \right\|^2} \le \frac{{({{\left\| {{x_N} - p} \right\|}^2} + \frac{a}{{1 - a}}{{[{{\left\| {{x_N} - p} \right\|}^2} - {{\left\| {{x_{N - 1}} - p} \right\|}^2}]}_ + } + \frac{{(2 + \frac{{1 - \theta }}{{4\theta }})M}}{{1 - a}})\frac{1}{{\theta (1 - \mu )}}}}{{n - N + 1}}.
\end{equation}
Since ${[\mathop {\min }\limits_{N \le i \le n} {\left\| {{w_i} - {y_i}} \right\|^2}]^{\frac{1}{2}}} = \mathop {\min }\limits_{N \le i \le n} \left\| {{w_i} - {y_i}} \right\|$, we get
\begin{equation}
	\label{+11}
	\mathop {\min }\limits_{N \le i \le n} \left\| {{w_i} - {y_i}} \right\| \le {(\frac{{({{\left\| {{x_N} - p} \right\|}^2} + \frac{a}{{1 - a}}{{[{{\left\| {{x_N} - p} \right\|}^2} - {{\left\| {{x_{N - 1}} - p} \right\|}^2}]}_ + } + \frac{{(2 + \frac{{1 - \theta }}{{4\theta }})M}}{{1 - a}})\frac{1}{{\theta (1 - \mu )}}}}{{n - N + 1}})^{\frac{1}{2}}}.
\end{equation}
This completes the proof.
\begin{remark}
	By Lemma \ref{l2.6} we know ${y_n} = {w_n}$ implies that $y_n$ is a solution of MIP. This means that the error estimate given in Theorem \ref{t3.2} can be regarded as the convergence rate of Algorithm \ref{alg3.1}.
\end{remark}
\section{Strong convergence}
\setcounter{equation}{0}
In this section, we analyse strong convergence and linear convergence of Algorithm \ref{alg3.1}. Firstly, we give the following assumption.
\begin{assumption}
	\label{a3.2}
The mapping $A:H \to H$ is $L$-Lipschitz continuous, $r$-strongly monotone and the set-valued mapping $B:H \to {2^H}$ is maximal monotone or the mapping $A:H \to H$ is $L$-Lipschitz continuous, monotone and the set-valued mapping $B:H \to {2^H}$ is $r$-strongly monotone.
\end{assumption}
\begin{theorem}	
	\label{t4.1}Assume that the conditions ($C_1$),($C_3$) and Assumption \ref{a3.2} hold. Let $\{ {x_n}\} $ be a sequence generated by Algorithm \ref{alg3.1}. Then $\{ {x_n}\} $ converges strongly to some solution $p \in \Omega $.
\end{theorem}	
	\begin{proof}
By \eqref{3.8} \eqref{3.9} and Lemma \ref{l3.2} (ii), we have
	\begin{equation}
		\label{3.36}
		\begin{array}{c}
			\begin{split}
				{\left\| {{x_{n + 1}} - p} \right\|^2} &\le (1 - {\theta _n}){\left\| {{z_n} - p} \right\|^2} + {\theta _n}{\left\| {{w_n} - p} \right\|^2} - {\theta _n}(1 - \frac{{{(\mu  + {\mu _n})^2}{\lambda _n}^2}}{{{\lambda _{n + 1}}^2}}){\left\| {{w_n} - {y_n}} \right\|^2}\\
				&\quad - 2{\theta _n}r{\lambda _n}{\left\| {{y_n} - p} \right\|^2}- \frac{{(1 - {\theta _n})}}{{{\theta _n}}}{\left\| {{x_{n + 1}} - {z_n}} \right\|^2}\\
				&\le (1 - {\theta _n}){\left\| {{z_n} - p} \right\|^2} + {\theta _n}{\left\| {{w_n} - p} \right\|^2} - {\theta _n}(1 - \frac{{{(\mu  + {\mu _n})^2}{\lambda _n}^2}}{{{\lambda _{n + 1}}^2}}){\left\| {{w_n} - {y_n}} \right\|^2}\\
				&\quad - 2{\theta _n}r{\lambda _n}{\left\| {{y_n} - p} \right\|^2}.
			\end{split}
		\end{array}
	\end{equation}
Lemma \ref{l3.1} and the fact$\mathop {\lim }\limits_{n \to \infty } {\mu _n} = 0$ imply that $\mathop {\lim }\limits_{n \to \infty } 1 - \frac{{{(\mu  + {\mu _n})^2}{\lambda _n}^2}}{{{\lambda _{n + 1}}^2}} = 1 - {\mu ^2} > 0$. Thus there exists a positive integer $N \ge 1$ such that
$${\theta _n}(1 - \frac{{{(\mu  + {\mu _n})^2}{\lambda _n}^2}}{{{\lambda _{n + 1}}^2}})\left\| {{w_n} - {y_n}} \right\|^2 \geq 0,~\forall~n \ge N .$$

It follows that from \eqref{3.36}
\begin{equation}
	\label{4.2}
	\begin{array}{c}
		\begin{split}
			{\left\| {{x_{n + 1}} - p} \right\|^2}&\le (1 - {\theta _n}){\left\| {{z_n} - p} \right\|^2} + {\theta _n}{\left\| {{w_n} - p} \right\|^2} - 2{\theta _n}r{\lambda _n}{\left\| {{y_n} - p} \right\|^2} \\
			&\le (1 - {\theta _n}){\left\| {{z_n} - p} \right\|^2} + {\theta _n}{\left\| {{w_n} - p} \right\|^2} - 2\theta r{\lambda _n}{\left\| {{y_n} - p} \right\|^2},~ \forall n \ge N.
		\end{split}
	\end{array}
\end{equation}
From Lemma \ref{l3.1}, it follows that
\begin{equation}
	\label{4.3}
	\begin{array}{c}
		{\left\| {{x_{n + 1}} - p} \right\|^2}\le (1 - {\theta _n}){\left\| {{z_n} - p} \right\|^2} + {\theta _n}{\left\| {{w_n} - p} \right\|^2} - 2\theta r\lambda^* {\left\| {{y_n} - p} \right\|^2},~ \forall n \ge N,
	\end{array}
\end{equation}
where $\lambda^*=\min \{ \frac{\mu }{L},{\lambda _1}\}$.\\
 The definition of $w_n$ implies that
\begin{equation}
	\label{4.4}
	\begin{array}{c}
		\begin{split}
			 	{\left\| {{w_n} - p} \right\|^2} &= {\left\| {{x_n} + {\alpha _n}({x_n} - {x_{n - 1}}) - p} \right\|^2}\\
			 	&\le {\left\| {{x_n} - p} \right\|^2} + {\alpha _n}^2{\left\| {{x_n} - {x_{n - 1}}} \right\|^2} + 2{\alpha _n}\left\| {{x_n} - p} \right\|\left\| {{x_n} - {x_{n - 1}}} \right\|.
		\end{split}
	\end{array}
\end{equation}
By utilizing the definition of $z_n$, we have
\begin{equation}
	\label{4.5}
	\begin{array}{c}
		{\left\| {{z_n} - p} \right\|^2} \le {\left\| {{x_n} - p} \right\|^2} + {\beta _n}^2{\left\| {{x_n} - {x_{n - 1}}} \right\|^2} + 2{\beta _n}\left\| {{x_n} - p} \right\|\left\| {{x_n} - {x_{n - 1}}} \right\|.
	\end{array}
\end{equation}
Combining  \eqref{4.3} \eqref{4.4} and \eqref{4.5}, we obtain
\begin{equation}
	\label{4.6}
	\begin{array}{c}
		\begin{split}
			 	{\left\| {{x_{n + 1}} - p} \right\|^2} &\le (1 - {\theta _n})({\left\| {{x_n} - p} \right\|^2} + {\beta _n}^2{\left\| {{x_n} - {x_{n - 1}}} \right\|^2} + 2{\beta _n}\left\| {{x_n} - p} \right\|\left\| {{x_n} - {x_{n - 1}}} \right\|)\\
			 	&\quad + {\theta _n}({\left\| {{x_n} - p} \right\|^2} + {\alpha _n}^2{\left\| {{x_n} - {x_{n - 1}}} \right\|^2} + 2{\alpha _n}\left\| {{x_n} - p} \right\|\left\| {{x_n} - {x_{n - 1}}} \right\|)\\
			 	&\quad - 2{\theta }\lambda^* r{\left\| {{y_n} - p} \right\|^2}\\
			 	&= {\left\| {{x_n} - p} \right\|^2} + ((1 - {\theta _n}){\beta_n}^2 + {\theta _n}{\alpha _n}^2){\left\| {{x_n} - {x_{n - 1}}} \right\|^2}\\
			 	&\quad + 2((1 - {\theta _n}){\beta _n} + {\theta _n}{\alpha _n})\left\| {{x_n} - p} \right\|\left\| {{x_n} - {x_{n - 1}}} \right\| - 2{\theta}\lambda^* r{\left\| {{y_n} - p} \right\|^2},~ \forall n \ge N.
		\end{split}
	\end{array}
\end{equation}
In addition,
\begin{equation*}
	\begin{array}{c}
		\begin{split}
			{\left\| {{x_n} - p} \right\|^2} &\le 2({\left\| {{x_n} - {y_n}} \right\|^2} + {\left\| {{y_n} - p} \right\|^2})\\
			&\le 4({\left\| {{x_n} - {w_n}} \right\|^2} + {\left\| {{y_n} - {w_n}} \right\|^2}) + 2{\left\| {{y_n} - p} \right\|^2},
		\end{split}
	\end{array}
\end{equation*}
which implies that
\begin{equation}
	\label{4.8}
	\begin{array}{c}
		\begin{split}
			{\left\| {{y_n} - p} \right\|^2} \ge \frac{1}{2}{\left\| {{x_n} - p} \right\|^2} - 2{\left\| {{y_n} - {w_n}} \right\|^2} - 2{\left\| {{w_n} - {x_n}} \right\|^2}.
		\end{split}
	\end{array}	
\end{equation}
In view of \eqref{4.6} and \eqref{4.8}, we have
\begin{equation}
	\label{4.9}
	\begin{array}{c}
		\begin{split}
				{\left\| {{x_{n + 1}} - p} \right\|^2} &\le {\left\| {{x_n} - p} \right\|^2} + ((1 - {\theta _n}){\beta _n}^2 + {\theta _n}{\alpha _n}^2){\left\| {{x_n} - {x_{n - 1}}} \right\|^2}\\
				&\quad + 2((1 - {\theta _n}){\beta _n} + {\theta _n}{\alpha _n})\left\| {{x_n} - p} \right\|\left\| {{x_n} - {x_{n - 1}}} \right\|\\
				&\quad- \theta r\lambda^* {\left\| {{x_n} - p} \right\|^2} + 4\theta r\lambda^* {\left\| {{y_n} - {w_n}} \right\|^2} + 4\theta r\lambda^* {\left\| {{w_n} - {x_n}} \right\|^2}\\
				&= (1 - \theta r\lambda^* ){\left\| {{x_n} - p} \right\|^2} + ((1 - {\theta _n}){\beta _n}^2 + {\theta _n}{\alpha _n}^2){\left\| {{x_n} - {x_{n - 1}}} \right\|^2}\\
				&\quad + 2((1 - {\theta _n}){\beta _n} + {\theta _n}{\alpha_n})\left\| {{x_n} - p} \right\|\left\| {{x_n} - {x_{n - 1}}} \right\| \\
				&\quad+ 4\theta r\lambda^* {\left\| {{y_n} - {w_n}} \right\|^2} + 4\theta r\lambda^* {\left\| {{w_n} - {x_n}} \right\|^2},~ \forall n \ge N.
		\end{split}
	\end{array}	
\end{equation}
Since ${\lambda ^*} \le \frac{\mu }{L}$ and $r\le L$, we obtain $1 - \theta r\lambda^*\in(0,1) $.

Since $x_n$ is bounded, and $\mathop {\lim }\limits_{n \to \infty } \left\| {{x_{n+1}} - {x_n}} \right\| = \mathop {\lim }\limits_{n \to \infty } \left\| {{x_n} - {w_n}} \right\| = \mathop {\lim }\limits_{n \to \infty } \left\| {{w_n} - {y_n}} \right\| = 0$, owing to Lemma \ref{l2.4}, we have
\begin{equation}
	\label{4.11}	
	\mathop {\lim }\limits_{n \to \infty } {\left\| {{x_n} - p} \right\|^2} =0.
\end{equation}
Hence
\begin{equation}
	\label{4.12}	
	\mathop {\lim }\limits_{n \to \infty } {\left\| {{x_n} - p} \right\|}= 0.
\end{equation}
This proof is completed.
\end{proof}
\begin{remark}
\label{r4.1}
To the best of our knowledge, Theorem \ref{t4.1} is one of the few available strong convergence results for algorithms with the double inertial extrapolation steps to solve MIP (\ref{1.1}). In addition, we emphasize that Theorem \ref{t4.1} does not need to know the modulus of strong monotonicity and the Lipschitz constant in advance.
\end{remark}
Owing to Proposition \ref{p3.1} and Theorem \ref{t4.1}, it is easy to get the following results.
\begin{corollary}
	\label{c4.1} Let $B=N_C$ and $A$ is $\mu$-strong pseudomonotone and $L$-Lipschitz continuous. If the conditions ($C_3$) hold, then the sequence $\{ {x_n}\} $ from Algorithm \ref{alg3.1} converges strongly to a point $p\in S$.
\end{corollary}
\begin{remark}
\label{r4.2}Corollary \ref{c4.1} improves Theorem 5.1 of \cite{YI} in the
following aspects: (i) the sequence $\{\alpha_n\}$ may not be non-decreasing; (ii) the sequence $\{\beta_n\}$ may not be a constant; (iii) we require $0<\theta<\theta_n\le\theta_{n+1}<\dfrac{1}{1+\epsilon},\varepsilon  \in (1, + \infty )$ other than $ \varepsilon  \in (2, + \infty )$, which extend the taking value interval of $\theta_n$.
\end{remark}
In order to discuss the linear convergence rate of our algorithm,  we need the following assumption.
\begin{assumption}
	\label{a4.1}
	\begin{itemize}
		\item[(i)] The solution set of the inclusion problem \eqref{1.1} is nonempty,  that is, $\Omega\ne \emptyset $.
		\item[(ii)] The mapping $A:H \to H$ is $L$-Lipschitz continuous, $r$-strongly monotone and the set-valued mapping $B:H \to {2^H}$ is maximal monotone or the mapping $A:H \to H$ is $L$-Lipschitz continuous, monotone and the set-valued mapping $B:H \to {2^H}$ is $r$-strongly monotone.
		\item[(iii)] Let $\hat{\lambda} :=\min \{ \frac{\mu }{L},{\lambda _1}\}$, $\tau  := 1 - \frac{1}{2}\min \{ 1 - \mu ,2\hat{\lambda} r \}  \in (\frac{1}{2},1)$ and the following conditions hold:
\end{itemize}
\begin{itemize}
		\item[($c_1$)] $0 \le {\beta _n} \le \beta  < \frac{1}{2}(\frac{1}{\tau } - 1) $
		\item[($c_2$)] $0 \le {\alpha_n}  \le \alpha< \frac{{1 - \tau }}{\tau }$
		\item[($c_3$)] $ \max\{\frac{1-\beta}{1+\alpha-\beta},~\frac{{\beta}  }{{1 + {\beta}   - \tau (1 + {\alpha} )}}\} <\theta \le  \theta_{n-1}\le  \theta_n  \le\frac{{ - 1 - \beta  + \sqrt {{{(1 + \beta )}^2} - 4(\frac{1}{\tau } - 1 - 2\beta )(\beta  - 1)} }}{{2(\frac{1}{\tau } - 1 - 2\beta )}}$
	\end{itemize}
\end{assumption}

\begin{remark}
	\label{r4.2}
	The parameters set satisfying the condition is non-empty. For example, we can choose $L=1.5$, $r=1$, $\mu=0.45$, $\beta=\beta_n=0.1$, $\alpha=\alpha_n=0.37$ and $\theta_n=0.72$.
\end{remark}
 Next, we establish the linear convergence rate of our algorithm under Assumption \ref{a4.1}.
\begin{theorem}	
	\label{t4.2} Suppose that Assumption \ref{a4.1} hold. Then $\{ {x_n}\} $ generated by Algorithm \ref{alg3.1} converges linearly to some point in $\Omega$.
\end{theorem}
\begin{proof}
By Lemma \ref{l3.1} we know  $\mathop {\lim }\limits_{n \to \infty } {\lambda _n} = \lambda \geq \hat{\lambda} =\min \{ \frac{\mu }{L},{\lambda _1}\}$. Since $\mathop {\lim }\limits_{n \to \infty } {\mu _n} =0$, $\mathop {\lim }\limits_{n \to \infty } 1 - \frac{{{(\mu  + {\mu _n})^2}{\lambda _n}^2}}{{{\lambda _{n + 1}}^2}} = 1 - {\mu ^2} >1 - \mu > 0$. Thus it follows that from \eqref{3.36} there exists a positive integer $N \ge 1$ such that
	\begin{equation}
	\label{4.14}
	\begin{array}{c}
		\begin{split}
			{\left\| {{x_{n + 1}} - p} \right\|^2}&\le (1 - {\theta_n  }){\left\| {{z_n} - p} \right\|^2} + {\theta_n  }{\left\| {{w_n} - p} \right\|^2}- \frac{{(1 - {\theta_n  })}}{{{\theta_n  }}}{\left\| {{x_{n + 1}} - {z_n}} \right\|^2}\\
			&\quad- {\theta_n  }(1 - {\mu }){\left\| {{w_n} - {y_n}} \right\|^2} - 2{\theta_n  }r\hat{\lambda}{\left\| {{y_n} - p} \right\|^2}\\
			&\le (1 - {\theta_n  }){\left\| {{z_n} - p} \right\|^2} + {\theta_n  }{\left\| {{w_n} - p} \right\|^2}- \frac{{(1 - {\theta_n  })}}{{{\theta_n }}}{\left\| {{x_{n + 1}} - {z_n}} \right\|^2}\\
			&\quad- \theta_n  \min \{ 1 - \mu ,2\hat{\lambda} r \} ( {\left\| {{w_n} - {y_n}} \right\|^2}+{\left\| {{y_n} - p} \right\|^2} ) \\
			&\le (1 - \theta_n ){\left\| {{z_n} - p} \right\|^2} + {\theta_n }\tau {\left\| {{w_n} - p} \right\|^2}- \frac{{(1 - {\theta_n  })}}{{{\theta_n }}}{\left\| {{x_{n + 1}} - {z_n}} \right\|^2},~\forall n \ge N.
		\end{split}
	\end{array}
\end{equation}
Substituting \eqref{3.13}, \eqref{3.14} and \eqref{3.15} into \eqref{4.14}, we have
	\begin{equation}
		\label{4.15}
		\begin{array}{c}
			\begin{split}
				{\left\| {{x_{n + 1}} - p} \right\|^2} &\le (1 - \theta_n )[(1 + {\beta_n}  ){\left\| {{x_n} - p} \right\|^2} - {\beta_n}  {\left\| {{x_{n - 1}} - p} \right\|^2}+ (1 + {\beta_n}  )\beta_n {\left\| {{x_n} - {x_{n - 1}}} \right\|^2}] \\
				&\quad+ \theta_n  \tau [(1 + {\alpha_n} ){\left\| {{x_n} - p} \right\|^2}- {\alpha_n}{\left\| {{x_{n - 1}} - p} \right\|^2} + (1 + {\alpha_n}){\alpha_n} {\left\| {{x_n} - {x_{n - 1}}} \right\|^2}]\\
				&\quad- \frac{{(1 - \theta_n  )}}{\theta_n  }[(1 - {\beta_n}  ){\left\| {{x_{n + 1}} - {x_n}} \right\|^2} + ({{\beta_n}  ^2} -{\beta_n}  ){\left\| {{x_n} - {x_{n - 1}}} \right\|^2}]\\
				&\le [(1 - \theta_n  )(1 + {\beta_n}  ) + \theta_n  \tau (1 + {\alpha_n} )]{\left\| {{x_n} - p} \right\|^2}- [(1 - \theta_n  ){\beta_n}   + \theta_n  \tau {\alpha_n}]{\left\| {{x_{n-1}} - p} \right\|^2} \\
				&\quad+ [(1 - \theta_n  )(1 + {\beta_n}  ){\beta_n}  + \theta_n  \tau (1 + {\alpha_n} ){\alpha_n}  + \frac{{(1 - \theta_n  )}}{\theta_n  }({\beta_n}   - {{\beta_n} ^2})]{\left\| {{x_n} - {x_{n - 1}}} \right\|^2}\\
				&\quad- \frac{{(1 - \theta_n  )}}{\theta_n  }(1 - {\beta_n} ){\left\| {{x_{n + 1}} - {x_n}} \right\|^2}\\
				&\le [(1 - \theta_n )(1 + {\beta_n}  ) + \theta_n \tau (1 + {\alpha_n} )]{\left\| {{x_n} - p} \right\|^2} \\
				&\quad+ [(1 - \theta_n  )(1 + {\beta_n}  ){\beta_n}  + \theta_n \tau (1 + {\alpha_n} ){\alpha_n}  + \frac{{(1 - \theta_n  )}}{\theta_n  }({\beta_n}   - {{\beta_n}  ^2})]{\left\| {{x_n} - {x_{n - 1}}} \right\|^2}\\
				&\quad- \frac{{(1 -\theta_n  )}}{\theta_n }(1 - {\beta_n}  ){\left\| {{x_{n + 1}} - {x_n}} \right\|^2}\\
				&\le [(1 - \theta_n )(1 + {\beta}  ) + \theta_n \tau (1 + {\alpha} )]{\left\| {{x_n} - p} \right\|^2} \\
				&\quad+ [(1 - \theta_n  )(1 + {\beta}  ){\beta}  + \theta_n \tau (1 + {\alpha} ){\alpha}  + \frac{{(1 - \theta_n  )}}{\theta_n  }({\beta}   - {{\beta}  ^2})]{\left\| {{x_n} - {x_{n - 1}}} \right\|^2}\\
				&\quad- \frac{{(1 -\theta_n  )}}{\theta_n }(1 - {\beta}  ){\left\| {{x_{n + 1}} - {x_n}} \right\|^2},~\forall n \ge N.
			\end{split}
		\end{array}
	\end{equation}
	This implies that
	\begin{equation}
		\label{4.16}
		\begin{array}{c}
				{\left\| {{x_{n + 1}} - p} \right\|^2} + \sigma_n{\left\| {{x_{n + 1}} - {x_n}} \right\|^2} \le [(1 - \theta_n  )(1 + {\beta}  ) + \theta_n \tau (1 + {\alpha})]({\left\| {{x_n} - p} \right\|^2} + \delta_n {\left\| {{x_n} - {x_{n - 1}}} \right\|^2}).
		\end{array}
	\end{equation}
	where $\delta_n  = \frac{{(1 - \theta_n  )(1 + {\beta}  ){\beta}  + \theta_n  \tau (1 + {\alpha}){\alpha}  + \frac{{(1 - \theta_n  )}}{\theta_n  }({\beta}   - {{\beta}  ^2})}}{{(1 - \theta_n  )(1 + {\beta}  ) + \theta_n  \tau (1 + {\alpha} )}}$ and $\sigma_n  = \frac{{(1 - \theta_n  )}}{\theta_n  }(1 - {\beta}  )$.\\
Now, we will show that $\delta_n  < \sigma_n $:
	\begin{equation}
		\label{4.17}
		\begin{array}{c}
			\begin{split}
				\delta_n  - \sigma_n  &= \frac{{(1 - \theta_n )(1 + {\beta}  ){\beta}  + \theta_n \tau (1 + {\alpha} ){\alpha}  + \frac{{1 - \theta_n }}{\theta_n }({\beta}   - {{\beta}  ^2})}}{{(1 - \theta_n )(1 + {\beta}  ) + \theta_n \tau (1 + {\alpha} )}} - \frac{{1 - \theta_n }}{\theta_n }(1 - {\beta}  )\\
						&= \frac{{(1 - {\theta _n})(1 + {\beta}){\beta} + {\theta _n}\tau (1 + {\alpha }){\alpha } + \frac{{1 - {\theta _n}}}{{{\theta _n}}}({\beta } - {\beta }^2)}}{{(1 - {\theta _n})(1 + {\beta }) + {\theta _n}\tau (1 + {\alpha })}}\\
						&\quad- \frac{{[(1 - {\theta _n})(1 + {\beta }) + {\theta _n}\tau (1 + {\alpha })]\frac{{1 - {\theta _n}}}{{{\theta _n}}}(1 - {\beta })}}{{(1 - {\theta _n})(1 + {\beta }) + {\theta _n}\tau (1 + {\alpha })}}\\
						&= \frac{{\frac{{1 - {\theta _n}}}{{{\theta _n}}}[{\theta _n}(1 + {\beta }){\beta } + {\beta } - {\beta }^2 - (1 - {\theta _n})(1 - {\beta }^2)]}}{{(1 - {\theta _n})(1 + {\beta }) + {\theta _n}\tau (1 + {\alpha })}}\\
						&\quad+ \frac{{\tau (1 + {\alpha })[{\theta _n}{\alpha } - (1 - {\theta _n})(1 - {\beta })]}}{{(1 - {\theta _n})(1 + {\beta }) + {\theta _n}\tau (1 + {\alpha })}}\\
						&= \frac{{\frac{{1 - {\theta _n}}}{{{\theta _n}}}[{\theta _n}{\beta } + {\beta } - 1 + {\theta _n}]}}{{(1 - {\theta _n})(1 + {\beta }) + {\theta _n}\tau (1 + {\alpha })}} + \frac{{\tau (1 + {\alpha })[{\theta _n}{\alpha } - 1 + {\theta _n} + {\beta } - {\theta _n}{\beta }]}}{{(1 - {\theta _n})(1 + {\beta }) + {\theta _n}\tau (1 + {\alpha })}}\\
						&= \frac{{[ - (1 + \beta ){\theta _n}^2 + 2{\theta _n} + \beta  - 1] + {\theta _n}\tau (1 + \alpha )[{\theta _n}(1 + \alpha ) - 1 + \beta  - {\theta _n}\beta ]}}{{{\theta _n}((1 - {\theta _n})(1 + \beta ) + {\theta _n}\tau (1 + \alpha ))}}\\
						&\le \frac{{[ - (1 + \beta ){\theta _n}^2 + 2{\theta _n} + \beta  - 1] + [\frac{1}{\tau }{\theta _n}^2 - {\theta _n} + \beta {\theta _n} - \beta {\theta _n}^2]}}{{{\theta _n}((1 - {\theta _n})(1 + \beta ) + {\theta _n}\tau (1 + \alpha ))}}\\
						&= \frac{{(\frac{1}{\tau } - 1 - 2\beta ){\theta _n}^2 + (1 + \beta ){\theta _n} + \beta  - 1}}{{{\theta _n}((1 - {\theta _n})(1 + \beta ) + {\theta _n}\tau (1 + \alpha ))}}
			\end{split}
		\end{array}
	\end{equation}
	Since $\theta_n  \le\frac{{ - 1 - \beta  + \sqrt {{{(1 + \beta )}^2} - 4(\frac{1}{\tau } - 1 - 2\beta )(\beta  - 1)} }}{{2(\frac{1}{\tau } - 1 - 2\beta )}}$, we get  $\delta_n  \leq \sigma_n $.  From $\theta_{n-1}\leq \theta_n$, it follows that $\sigma_n\leq\sigma_{n-1}$. By \eqref{4.16}, we have
	\begin{equation}
		\label{4.18}
		\begin{array}{c}
			\begin{split}
				{\left\| {{x_{n + 1}} - p} \right\|^2} + \sigma_n{\left\| {{x_{n + 1}} - {x_n}} \right\|^2} &\le [(1 - \theta_n )(1 + {\beta}) + \theta_n \tau (1 + {\alpha} )]({\left\| {{x_n} - p} \right\|^2} + \sigma_n {\left\| {{x_n} - {x_{n - 1}}} \right\|^2})\\
				&= [(1 +{\beta}) + \theta_n[ \tau (1 + {\alpha} )-(1+\beta)]({\left\| {{x_n} - p} \right\|^2} + \sigma_n {\left\| {{x_n} - {x_{n - 1}}} \right\|^2}).\\
&< [(1 +{\beta}) + \theta[ \tau (1 + {\alpha} )-(1+\beta)]({\left\| {{x_n} - p} \right\|^2} + \sigma_{n-1} {\left\| {{x_n} - {x_{n - 1}}} \right\|^2}),\\
				\end{split}
		\end{array}
	\end{equation}
where the last inequality follows from $\sigma_n\leq\sigma_{n-1}$ and $\alpha< \frac{{1 - \tau }}{\tau } - \frac{\beta }{{\tau }}$.\\
Since $\frac{\beta }{{(1 + \beta ) - \tau (1 + \alpha )}} < \theta $, we have $(1 +{\beta}) + \theta[ \tau (1 + {\alpha} )-(1+\beta)] \in (0,1)$. Therefore we get
	\begin{equation}
		\label{4.19}
		\begin{array}{c}
			{\left\| {{x_{n + 1}} - p} \right\|^2} + \sigma_n{\left\| {{x_{n + 1}} - {x_n}} \right\|^2}
			\le [(1 +{\beta}) + \theta[ \tau (1 + {\alpha} )-(1+\beta)]^{n-1}({\left\| {{x_2} - p} \right\|^2} + \sigma_1 {\left\| {{x_2} - {x_{1}}} \right\|^2}).
		\end{array}
	\end{equation}
This implies that
	\begin{equation}
		\label{4.20}
		{\left\| {{x_{n + 1}} - p} \right\|^2} \le [(1 +{\beta}) + \theta[ \tau (1 + {\alpha} )-(1+\beta)]^{n-1}({\left\| {{x_2} - p} \right\|^2} + \sigma_1 {\left\| {{x_2} - {x_{1}}} \right\|^2}).
	\end{equation}
Hence, the proof is completed.
\end{proof}	
\begin{remark}
Theorem 6.2 of \cite{YI} gives a linear convergence rate result of the subgradient extragradient method with double inertial
steps for solving variational inequalities. However, it only discuss the single inertial case in \cite{YI}.  As far as we know, algorithms with the double inertial extrapolation steps for solving variational inequalities and monotone inclusions have no linear convergent result. Hence Theorem \ref{4.2} is a new result.
\end{remark}

\section{Numerical experiments}
\setcounter{equation}{0}

\quad \quad In this section, we provide some numerical examples to show the performance of our Algorithm 3.1 (shortly Alg1), and compare it with others, including Abubakar et al's Algorithm 3.1 (shortly AKHAlg1) \cite{AB}, Chlolamjiak et al's Algorithm 3.2(CHCAlg2) \cite{CV}, Cholamjiak et al's Algorithm 3.1 (CHMAlg1) \cite{CH}, Hieu et al's Algorithm 3.1 (HAMAlg1) \cite{VA} and Yao et al's Algorithm 1 (YISAlg1) \cite{YI}.

All the programs were implemented in MATLAB R2021b on Intel(R) Core(TM) i7-7700HQ CPU@\\
2.80GHZ computer with RAM 8.00GB. We denote the number of iterations by "Iter." and the CPU time seconds by "CPU(s)".

\begin{example} 
	\label{E5.1}
We consider the signal recovery in compress sensing. This problem can be modeled as
\begin{equation}
	\label{5.1}
	y=Ax+\varepsilon.
\end{equation}
where ${\rm{y}} \in {{\rm{R}}^M}$ is observed or measured data, $A:{R^N} \to {R^M}$ is bounded linear operator, $x \in {R^N}$ is a vector with $K$ $(K<<N)$ nonzero components and $\varepsilon$ is the noise. It is know that the problem \eqref{5.1} can be viewed as the LASSO problem \cite{CV}
\begin{equation}
	\label{5.2}
	\mathop {\min }\limits_{x \in {R^N}} \frac{1}{2}\left\| {y - Ax} \right\|_2^2 + \lambda {\left\| x \right\|_1}~(\lambda>0).
\end{equation}
 The minimization problem \ref{5.2} is equivalent to the following monotone inclusion problem
\begin{equation}
	\label{5.3}
	{\rm find}~x \in {R^N}~{\rm such}~{\rm that}~0 \in (B + C)x.
\end{equation}
where $B = {A^T}(Ax - y)$ and $C = \partial (\lambda {\left\| x \right\|_1})$. In this experiment, the vector $x \in {R^N}$ is from uniform distribution in the interval $[-1,1]$.

The matrix $A \in {R^{M \times N}}$ is produced by a normal distribution with mean zero and one variance. The vector $y$ is generated by Gaussian noise $\varepsilon$ with variance $10^{-4}$. The initial points ${x_0}$ and ${x_1}$ are both zero. We use ${E_n} = \left\| {{x_n} - {x_{n - 1}}} \right\|$ to measure the restoration accuracy. And the stopping criterion is ${E_n} \le {10^{ - 5}}$.

 In the first experiment we consider the influence of different $\alpha_n$ and $\beta_n$ on the performance of our algorithm. We take $\mu  = 0.9$, ${\theta _n} = 0.45$, ${\lambda _1} = 0.1$, ${\mu _n} = 0$ and $p_n = \frac{1}{{{n^2}}}$.

 \begin{table}
  \caption{The performances of Algorithm 3.1 for different values of $\alpha_n$ and $\beta_n$ in Example 5.1}
  \centering
  \begin{tabular}{|l|l|l|l|l|l|l|}
   \hline
   \diagbox{$\beta_n$}{Iter.}{$\alpha_n$} & 0.2& 0.4&0.6&0.8&0.9&1\\
   \hline
   0&966&872&777&681&632&584\\
   \hline
   0.02&954&859&764&668&620&572\\
   \hline
   0.04&942&848&752&656&608&559\\
   \hline
   0.06&930&836&740&644&596&547\\
   \hline
   0.08&918&823&728&632&583&535\\
   \hline
   0.1&906&811&716&620&571&522\\
   \bottomrule
  \end{tabular}
 \end{table}
It can be seen from Table 1 that Algorithm 3.1 with double inertial extrapolation steps, that is,  $\beta_n\neq 0$  outperform the one with single inertia,
moreover, the increase of $\alpha_n$ and $\beta_n$ significantly improves the convergence
speed of the algorithm. This implies that it is important to investigate double inertial methods from both theoretical and numerical viewpoint.

 In second experiment we compared the performance of our algorithm with other algorithms. The following two cases are considered
\par\textbf{Case 1:} $K=20,~M=256,~N=512$;
\par\textbf{Case 2:} $K=40,~M=512,~N=1024$.

The parameters for algorithms are chosen as
	\par Alg1: $\mu  = 0.9$,${\alpha _n} = 1 - \frac{1}{{{{10}^n}}}$, ${\beta _n} = 0.1 - \frac{1}{{1000 + n}}$, ${\theta _n} = 0.45 - \frac{1}{{1000 + n}}$, ${\lambda _1} = 0.1$, ${\mu _n} = \frac{1}{{{n^2}}}$ and $p_n = \frac{1}{{{n^2}}};$
	\par CHCAlg2: $\mu  = 0.9$, $\alpha  = 0.1$, $\theta  = 1$ and ${\lambda _1} = 1;$
	\par CHMAlg1: $\mu  = 0.4$, ${\alpha _1} = 0.01$, ${\alpha _2} = 0.02$ and ${\lambda _0} = 0.01$;
	\par HAMAlg1: $\mu  = 0.4$, ${\lambda _{ - 1}} = {\lambda _0} = 0.1$;
	\par AKHAlg1: $\mu  = 0.3$, $\rho=0.1$, $\varrho=0.9$ and ${\lambda _0} = 1$.

Fig.1 and Table.2 show the numerical results of our algorithm and other algorithms in two cases respectively. We give the graphs of original signal and recovered signal in Fig. 2.

\begin{table}
	\caption{Numerical results for Example \ref{E5.1} }
	\centering
	\begin{tabular}{ccccc}
		\toprule
		\multirow{2}*{Algorithms} & Case 1& &Case 2&  \\
		\cline{2-5}
		& Iter. &CPU(s)&Iter.  &CPU(s)\\
		\midrule
		Alg1& 525&0.2447&809&1.2297\\
		CHCAlg2& 1347&0.5340&2595&3.0049\\
		CHMAlg1& 887&0.3496&1554&1.9651\\
		HAMAlg1& 971&0.3904&1577&1.5138\\
		AKHAlg1& 3427&1.3142&6723&6.4723\\
		\bottomrule
	\end{tabular}
\end{table}

\begin{figure}
	\centering
	\subfigure[$K=20$, $M=256$, $N=512$]{\includegraphics[width=7.5cm]{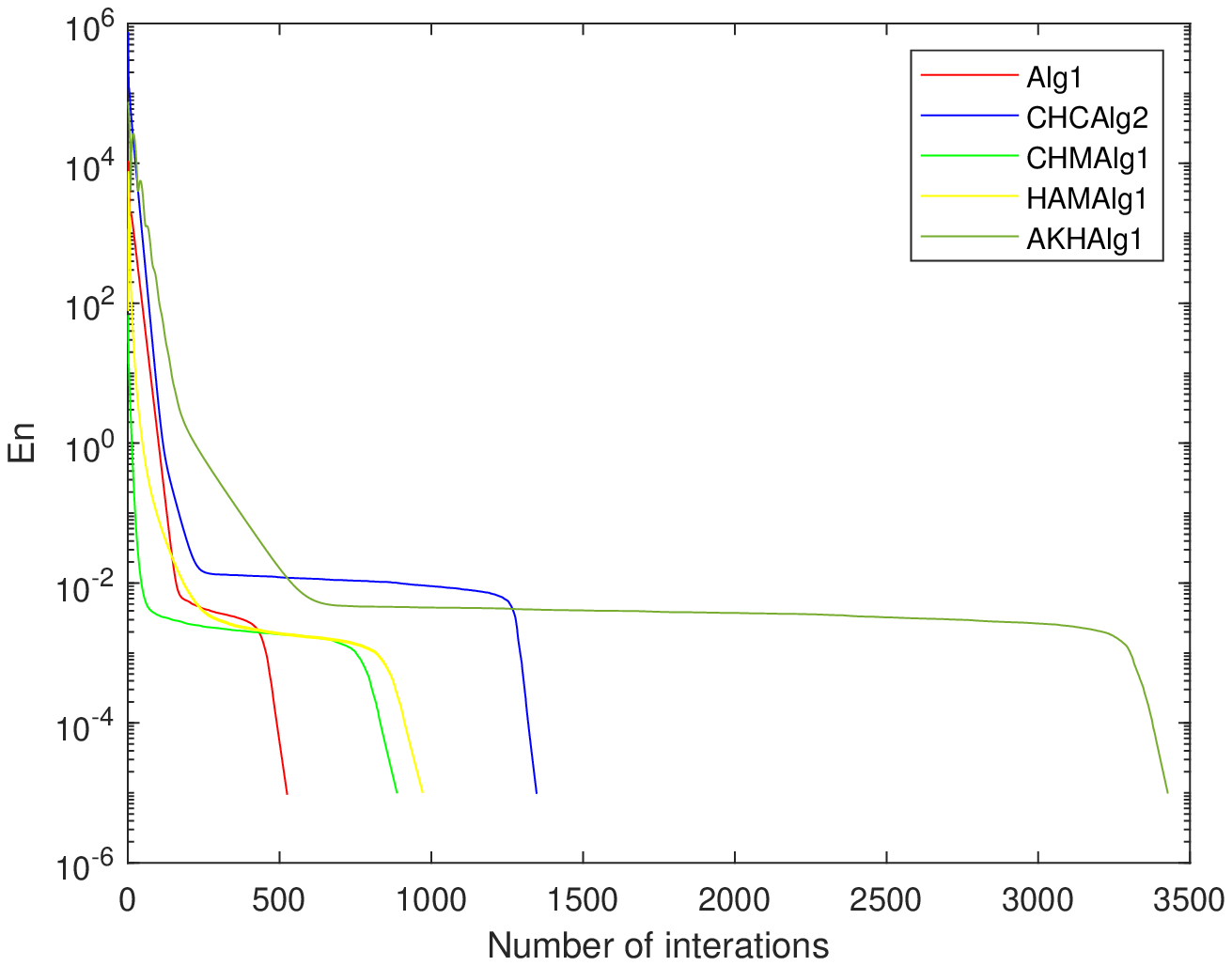}}
	\subfigure[$K=40$, $M=512$, $N=1024$]{\includegraphics[width=7.5cm]{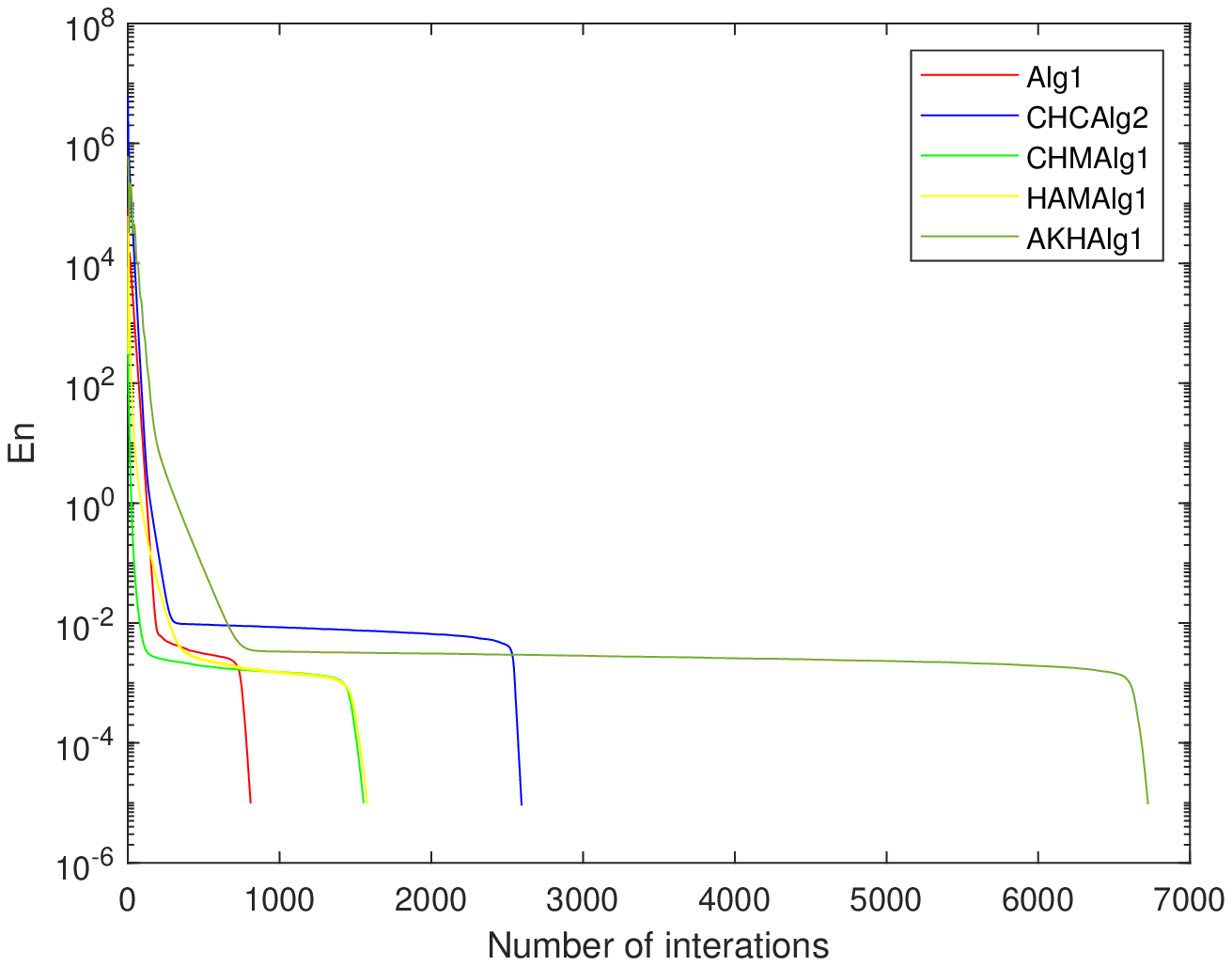}}
	\caption{Numerical behavior of $E_n$ for Example \ref{E5.1}}
\end{figure}
\begin{figure}
	\centering
	\subfigure[Original Signal(N=512, M=256, 20 spikes)]{\includegraphics[width=15cm]{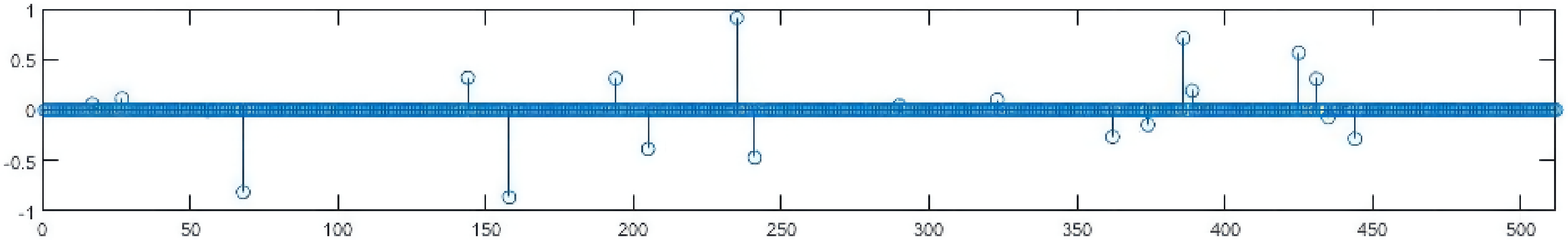}}
	\\ 
	\centering
	\subfigure[Measured values with variance $10^{-4}$]{\includegraphics[width=14.8cm]{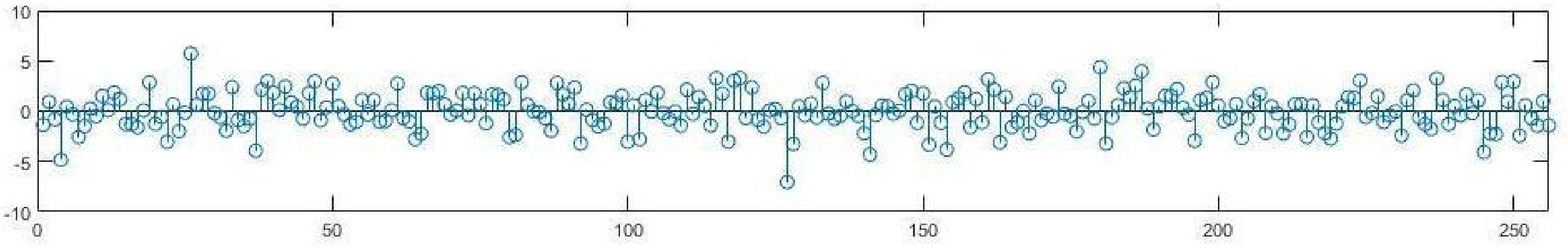}}
	\\ 
	\centering
	\subfigure[Recovered signal by method (Alg1)(484 iterations) ]{\includegraphics[width=15cm]{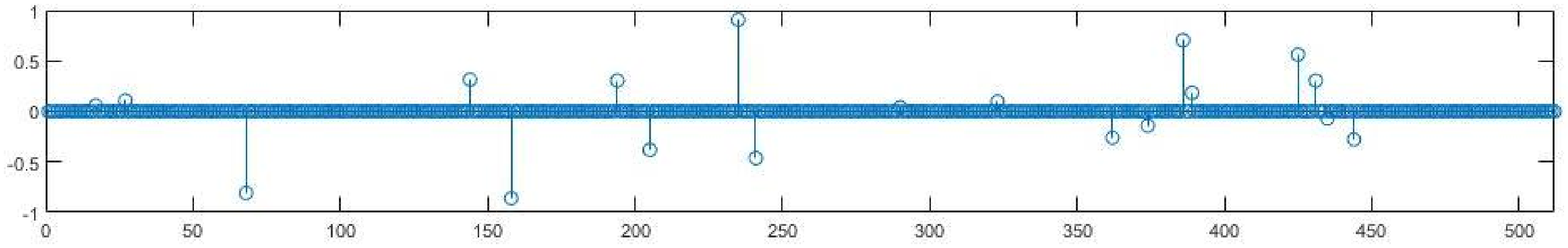}}
	\\ 
	\centering
	\subfigure[Recovered signal by method (CHCAlg2)(1347 iterations) ]{\includegraphics[width=15cm]{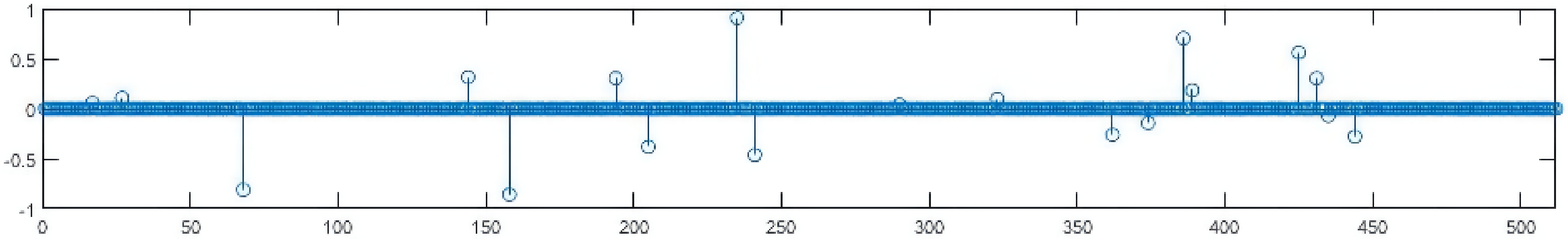}}
	\\ 
	\centering
	\subfigure[Recovered signal by method (CHMAlg1)(718 iterations)]{\includegraphics[width=15cm]{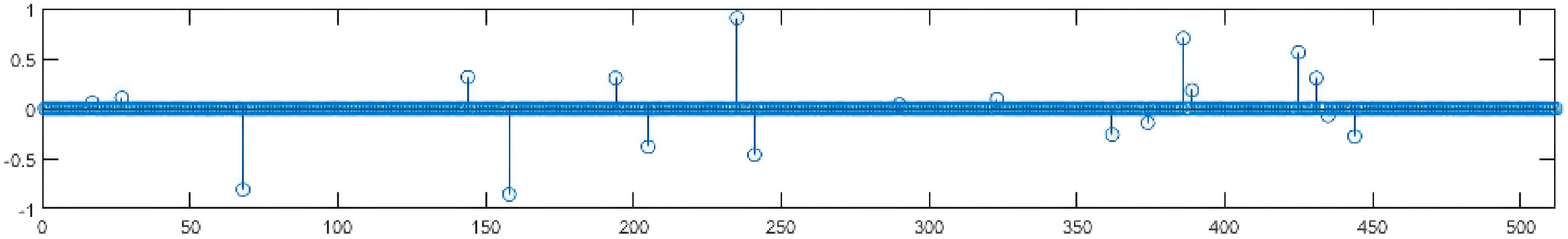}}
	\\ 
	\centering
	\subfigure[Recovered signal by method (HAMAlg1)(971 iterations)]{\includegraphics[width=15cm]{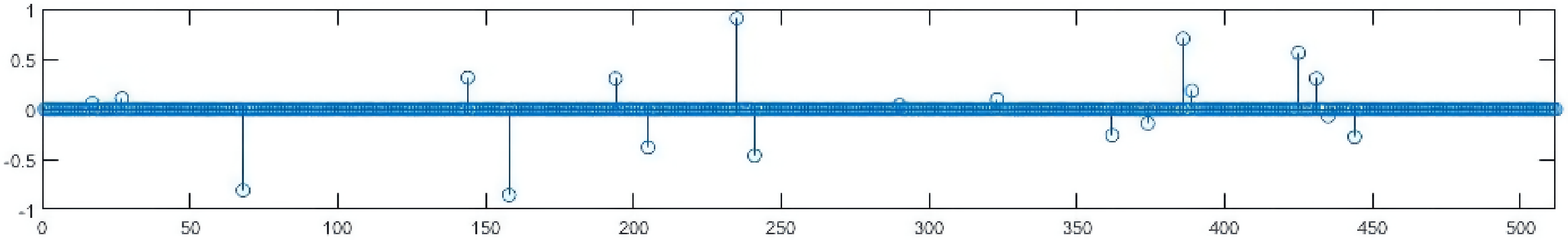}}
	\\ 
	\centering
	\subfigure[Recovered signal by method (AKHAlg1)(3427 iterations)]{\includegraphics[width=15cm]{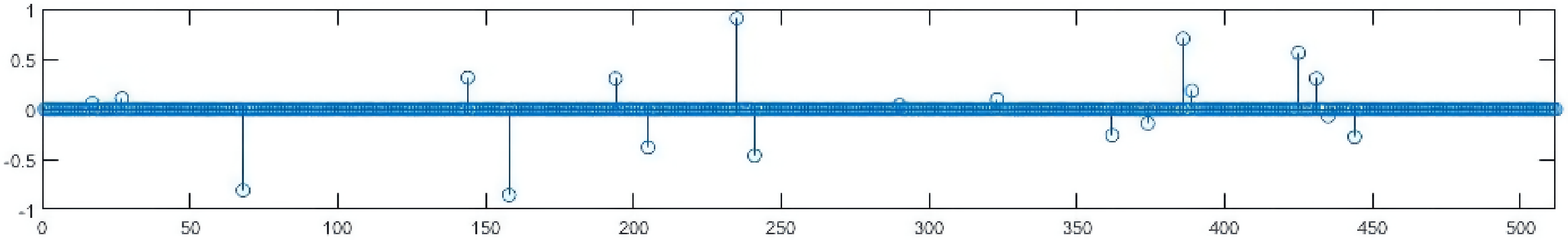}}
	\caption{From top to bottom: original signal, observation data, recovered signal by the methods (Alg1), (CHCAlg2), (CHMAlg1), (HAMAlg1) and (AKHAlg1) in Case 1, respectively}
\end{figure}
\begin{remark}
	\label{r5.1} It can be seen from Fig.1 and Table 2 that the number of iterations and CPU time of our algorithm are better than other algorithms, indicating that our algorithm has better performance.
\end{remark}

\end{example}

\begin{example} 
	\label{E5.2} We consider our algorithm to solve the variational inequality problem . The operator $A:{R^m} \to {R^m}$ is defined by $A(x): = Mx + q$ with  $M = N{N^T} + S + D$  and $q \in {R^m}$, where $N$ is an $m \times m$ matrix, $S$ is an $m \times m$ skew-symmetric matrix and $D$ is an $m \times m$ diagonal matrix. All entries of $N$ and $S$ are uniformly generated from $(-5,5)$ and all diagonal entries of $D$ are uniformly generated from $(0,0.3)$. It is easy to see that $M$ is positive definite. Define the feasible set $C:={R_m}^ + $ and use ${E_n} = \left\| {{x_n}} \right\|$ to measure the accuracy, and the stoping criterion is ${E_n} \le 10^{ - 3}$.
\par In the first experiment we consider the effect of relaxation coefficients $\theta_n$  on the performance of Algorithm 3.1.
	 We choose $\mu  = 0.9$, ${\alpha _n} = 1$, ${\beta _n} = 0.1$, ${\lambda _1} = 0.1$, ${\mu _n} = 0$ and $p_n = \frac{1}{{{n^2}}}$.
	
	\begin{table}
		\caption{The performances of Algorithm 3.1 for different values of $\theta_n$ in Example 5.2 }
		\centering
		\begin{tabular}{|l|l|l|l|l|l|l|l|l|l|}
			\hline
			 $\theta_n$&0.05&0.1&0.15&0.2&0.25&0.3&0.35&0.4&0.45\\
			\hline
			$Iter.$&16988&8521&5562&4035&3095&2454&1987&1360&1346\\
			\hline
			$CPU(s)$&1.8752&0.9888&0.6785&0.4804&0.3879&0.2910&0.2436&0.2246&0.1674\\
			\bottomrule
		\end{tabular}
	\end{table}
	Table 3 shows that the performance comparison of Algorithm 3.1 for different values of $\theta_n$. It can be seen that the performance of the algorithm becomes better with the increase of the relaxed parameter $\theta_n$. This indicates that the increase of the relaxation coefficient value range is of great significance to improve the performance of the algorithm.
	
	In second experiment we compared the performance of our algorithm with other algorithms. The following cases are considered
	\par\textbf{Case 1:} $m=50$; \quad\textbf{Case 2:} $m=100$; \quad\textbf{Case 3:} $m=150$; \quad\textbf{Case 4:} $m=200.$
	
	The parameters for algorithms are chosen as
	\par Alg1: $\mu  = 0.9$, ${\alpha _n} = 1 - \frac{1}{{{{10}^n}}}$, ${\beta _n} = 0.1 - \frac{1}{{1000 + n}}$, ${\theta _n} = 0.45 - \frac{1}{{1000 + n}}$, ${\mu _n} = 0$ and ${p _n} = \frac{1}{{{n^2}}};$
	\par HAMAlg1: $\mu  = 0.4$, ${\lambda _{ - 1}} = {\lambda _0} = 0.3;$
	\par CHCAlg2: $\mu  = 0.9$, $\alpha  = 0.3$, $\theta  = 0.4$ and ${\lambda _1} = 1;$
	\par YISAlg1: $\mu  = 0.9$, ${\alpha _n} = 0.2903$, $\delta  = 0.0241$, ${\theta _n} = 1$ and ${\lambda _1} = 0.1;$
	\par AKHAlg1: $\mu  = 0.3$, $\rho=0.1$, $\varrho=0.9$ and ${\lambda _0} = 1$;
	\par CHMAlg1: $\mu  = 0.4$, ${\alpha _1} = 0.01$, ${\alpha _2} = 0.02$ and ${\lambda _0} = 1$.
	
Fig.3 and Table.4 show the performance comparison of our algorithm with other algorithms in four cases respectively.
	\begin{table}
		\caption{Numerical results for Example \ref{E5.2}}
		\centering
		\begin{tabular}{ccccccccc}
			\toprule  
			\multirow{2}*{Algorithms}& $m=50$& &  $m=100$& &  $m=150$& &  $m=200$& \\
			\cline{2-9}
			& Iter. &CPU(s)&Iter.  &CPU(s)&Iter.  &CPU(s)&Iter. & CPU(s)\\
			\midrule  
			Alg1& 448&0.0191&642&0.0380&759&0.0968&1012&0.1681\\
			CHCAlg2& 723&0.0224&1048&0.0499&1234&0.1105&1644&0.2128\\
			YISAlg1& 600&0.0917&963&0.0635&1214&0.1374&1595&0.2125\\
		    HAMAlg1& 779&0.0203&1142&0.0475&1355&0.1410&1751&0.1940\\
		    AKHAlg1& 1049&0.0267&1425&0.0659&1691&0.1625&2371&0.2839\\
		    CHMAlg1& 871&0.0331&1182&0.0637&1408&0.1366&1879&0.2158\\
			\bottomrule

		\end{tabular}
	\end{table}
\begin{figure}
	\centering
	\subfigure[Case 1 $m=50$]{\includegraphics[width=7cm]{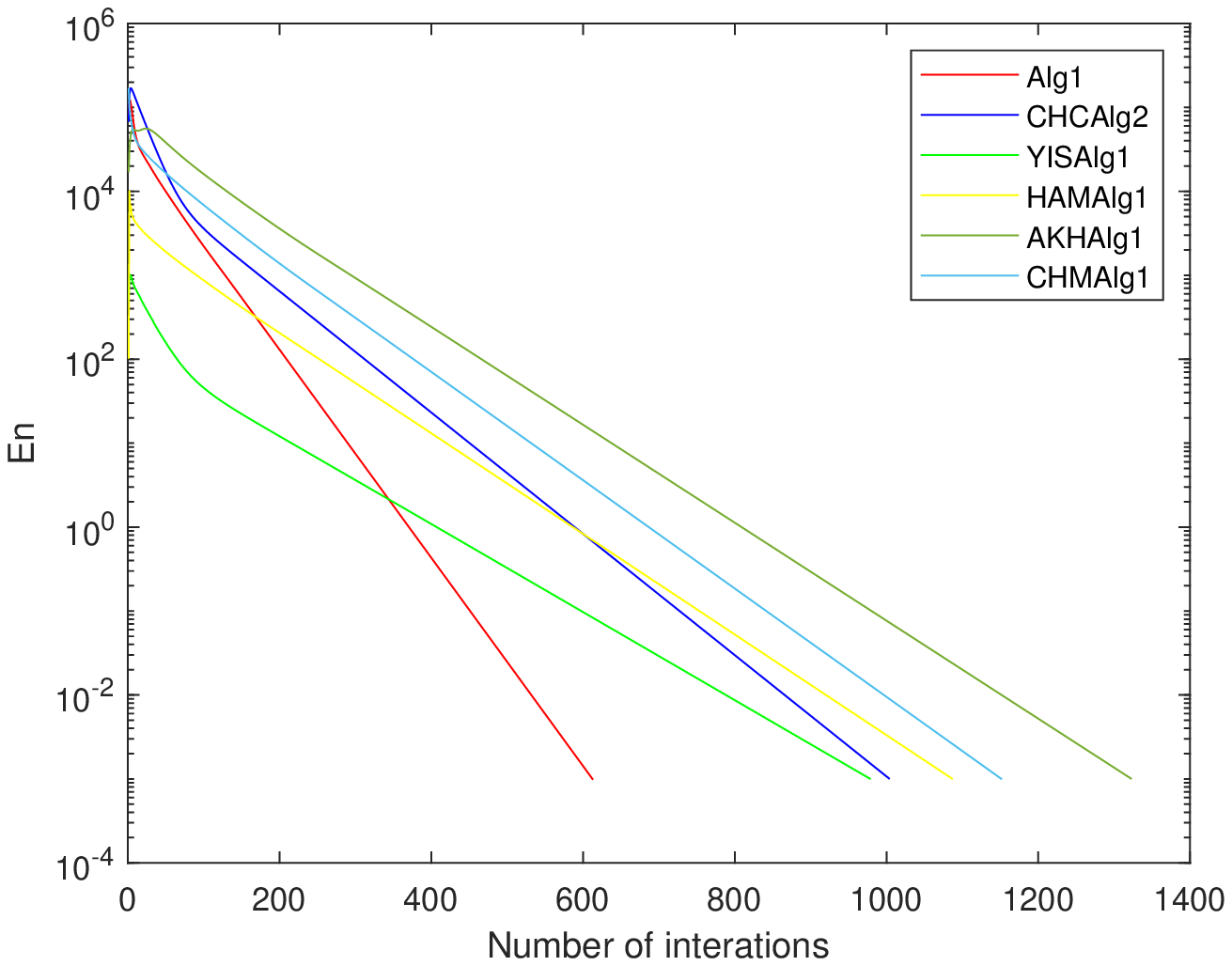}}
	\subfigure[Case 2 $m=100$]{\includegraphics[width=7cm]{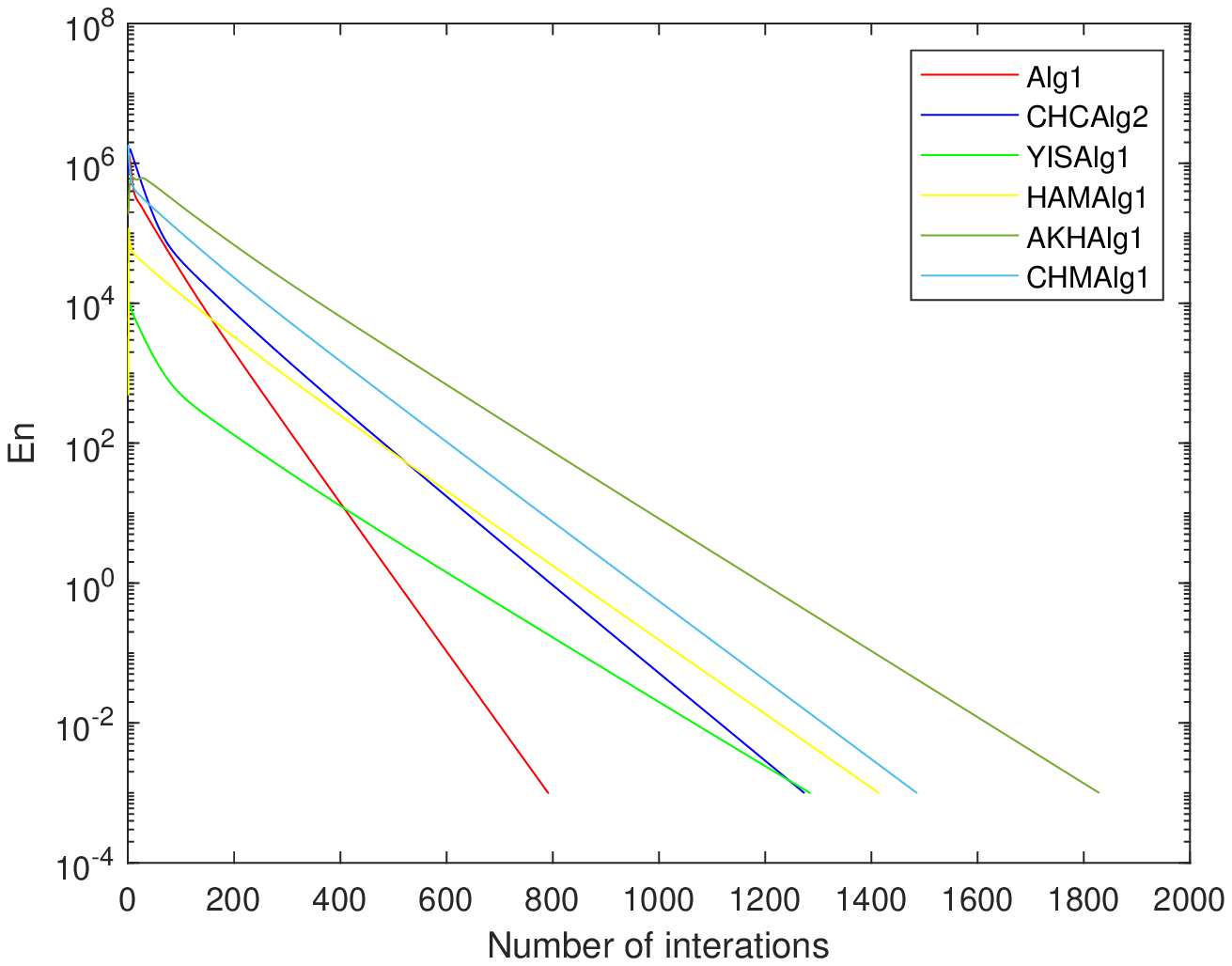}}
	\\ 
	\centering
	\subfigure[Case 3 $m=150$]{\includegraphics[width=7cm]{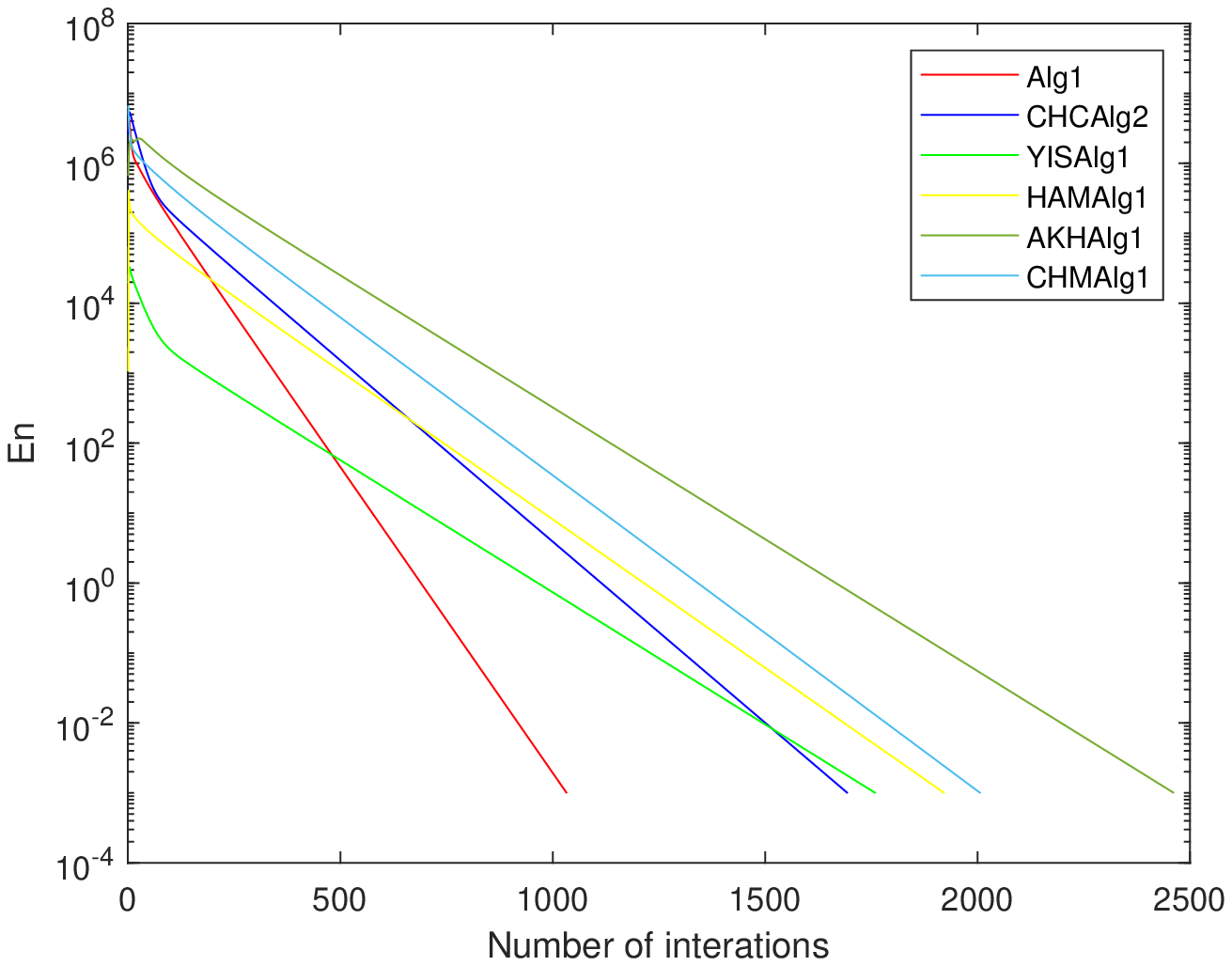}}
	\subfigure[Case 4 $m=200$]{\includegraphics[width=7cm]{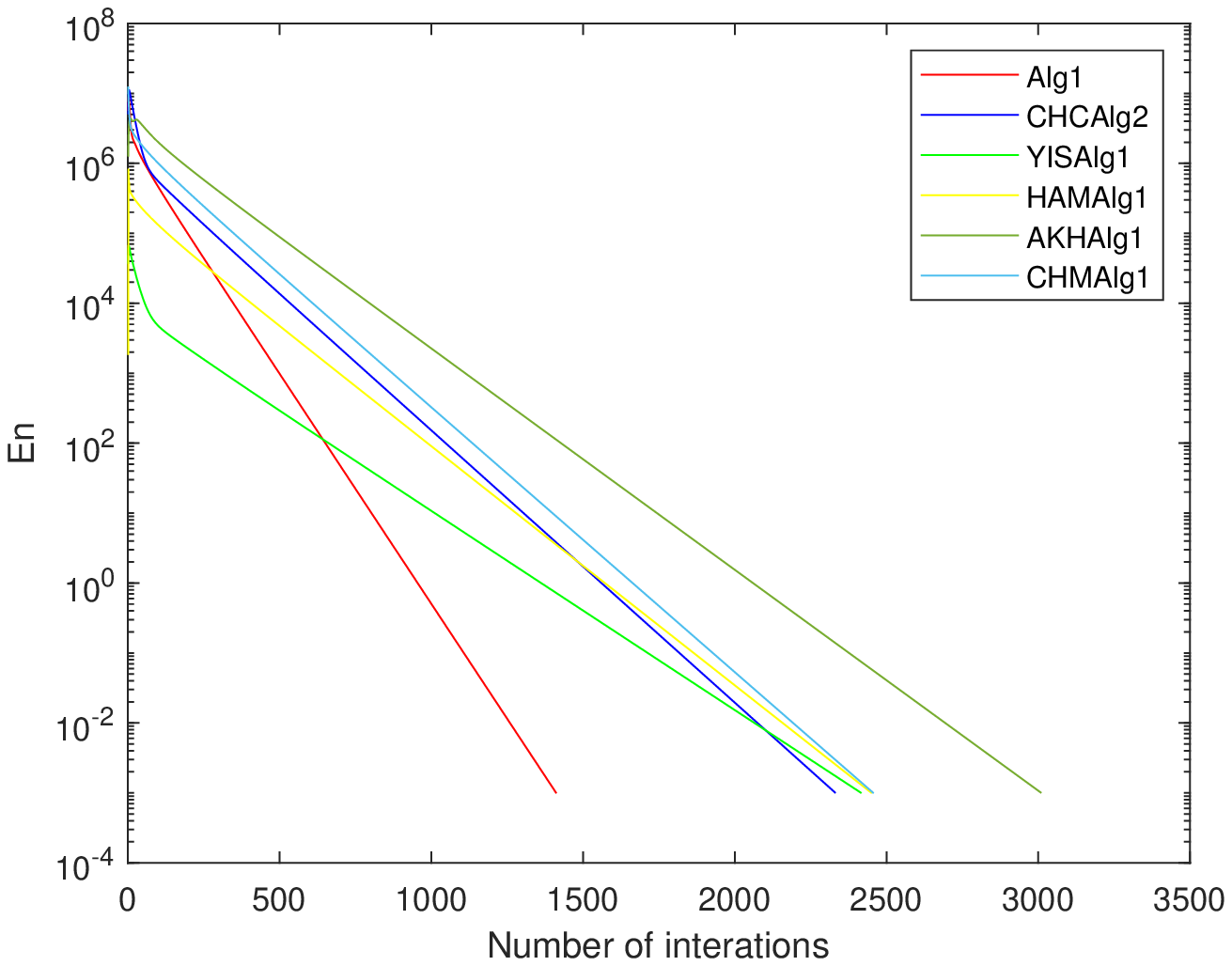}}
	\caption{Numerical behavior of $E_n$ for Example \ref{E5.2}} 
\end{figure}	
	\begin{remark}
		\label{r5.2} It can be seen that our algorithm performs better than other algorithms for such problems as Example 5.2 from Fig.3 and Table 4.
	\end{remark}
\end{example}

\begin{example} 
	\label{E5.3} Let $H: = {L_2}([1,2])$ with the norm
	\begin{equation*}
		\label{5.4}
			\left\| x \right\|: = {(\int_0^1 {x(t} )^2}dt{)^{\frac{1}{2}}}.
	\end{equation*}
	and the inner product
	\begin{equation*}
	\label{5.5}
	\left\langle {x,y} \right\rangle : = \int_0^1 {x(t)y(t)dt}.
\end{equation*}	\\
Let $C: = \{ x \in {L^2}([0,1]):\int_0^1 {tx(t)dt = 2} \} $ and  define $A:{L^2}([0,1]) \to {L^2}([0,1])$ by
	\begin{equation*}
	\label{5.5}
    Ax(t): = \max \{ x(t),0\} ,x \in {L^2}([0,1]),t \in [0,1].
\end{equation*}	
It is clear that $A$ is monotone and Lipschitz with $L=1$. The orthogonal projection onto $C$ have the following explicit formula
	\begin{equation*}
	\label{5.5}
	{P_C}(x)(t): = x(t) - \frac{{\int_0^1 {tx(t)dt - 2} }}{{\int_0^1 {{t^2}dt} }}t.
\end{equation*}	
The example is taken from \cite{YI}. Use ${E_n} = \left\| {{x_{n + 1}} - {x_n}} \right\|$ to measure the accuracy. The stoping criterion is ${E_n} \leq 10^{-4}$. \par The following cases are considered
\par\textbf{Case 1:} ${x_0} = \frac{{97{t^2} + 4t}}{{13}},{x_1} = \frac{{{t^2} - {e^{ - 7t}}}}{{250}};$
\par\textbf{Case 2:} ${x_0} = \frac{{97{t^2} + 4t}}{{13}},{x_1} = \frac{{\sin (3t) + \cos (10t)}}{{100}};$
\par\textbf{Case 3:} ${x_0} = \frac{{{t^2} - {e^{ - 7t}}}}{{250}},{x_1} = \frac{{\sin (3t) + \cos (10t)}}{{100}};$
\par\textbf{Case 4:} ${x_0} = \frac{{\sin (3t) + \cos (10t)}}{{100}},{x_1} = \frac{{97{t^2} + 4t}}{{13}}.$

The parameters for algorithms are chosen as
\par Alg1: $\mu  = 0.4$, ${\alpha _n} = 1 - \frac{1}{{{{10}^n}}}$, ${\beta _n} = 0.1 - \frac{1}{{1000 + n}}$, ${\theta _n} = 0.45 - \frac{1}{{1000 + n}}$ and ${p _n} = \frac{1}{{{n^2}}};$
\par CHCAlg2: $\mu  = 0.4$, $\alpha  = 0.3$, $\theta  = 0.4$ and ${\lambda _1} = 1;$
\par HAMAlg1: $\mu  = 0.4$, ${\lambda _{ - 1}} = {\lambda _0} = 0.1;$
\par YISAlg1: $\mu  = 0.4$, ${\alpha _n} = 0.2250$, $\delta  = 0.4950$, ${\theta _n} = 1$ and ${\lambda _1} = 1.1;$
\par AKHAlg1: $\mu  = 0.4$,$\rho=0.45$, $\varrho=0.3$ and ${\lambda _0} = 0.5;$
\par CHMAlg1: $\mu  = 0.4$, ${\alpha _1} = 0.01$, ${\alpha _2} = 0.02$ and ${\lambda _0} = 1$.

\begin{table}
	\caption{Numerical results for Example \ref{E5.3}}
	\centering
	\begin{tabular}{ccccccccc}
		\toprule  
		\multirow{2}*{Algorithms}& Case1& & Case2& &  Case3& & Case4& \\
		\cline{2-9}
		& Iter. &CPU(s)&Iter.  &CPU(s)&Iter.  &CPU(s)&Iter. & CPU(s)\\
		\midrule  
		Alg1   & 32 &0.0554 &32 &0.0499 &18 &0.0296 &36 &0.0416\\
		CHCAlg2& 40 &0.0601 &40 &0.0616 &24 &0.0361 &52 &0.0664\\
		VAMAlg1& 47 &0.0616 &46 &0.0558 &22 &0.0311 &70 &0.0974\\
		YISAlg1& 45 &0.0894 &45 &0.0629 &26 &0.0457 &53 &0.0912\\
		AKHAlg1& 38 &0.0746 &39 &0.0506 &20 &0.0277 &46 &0.0697\\
		HAMAlg1& 34 &0.0651 &37 &0.0498 &24 &0.0343 &70 &0.0802\\
		\bottomrule

	\end{tabular}
\end{table}

	\begin{figure}
		\centering
		\subfigure[${x_0} = \frac{{97{t^2} + 4t}}{{13}},{x_1} = \frac{{{t^2} - {e^{ - 7t}}}}{{250}}.$]{\includegraphics[width=7cm]{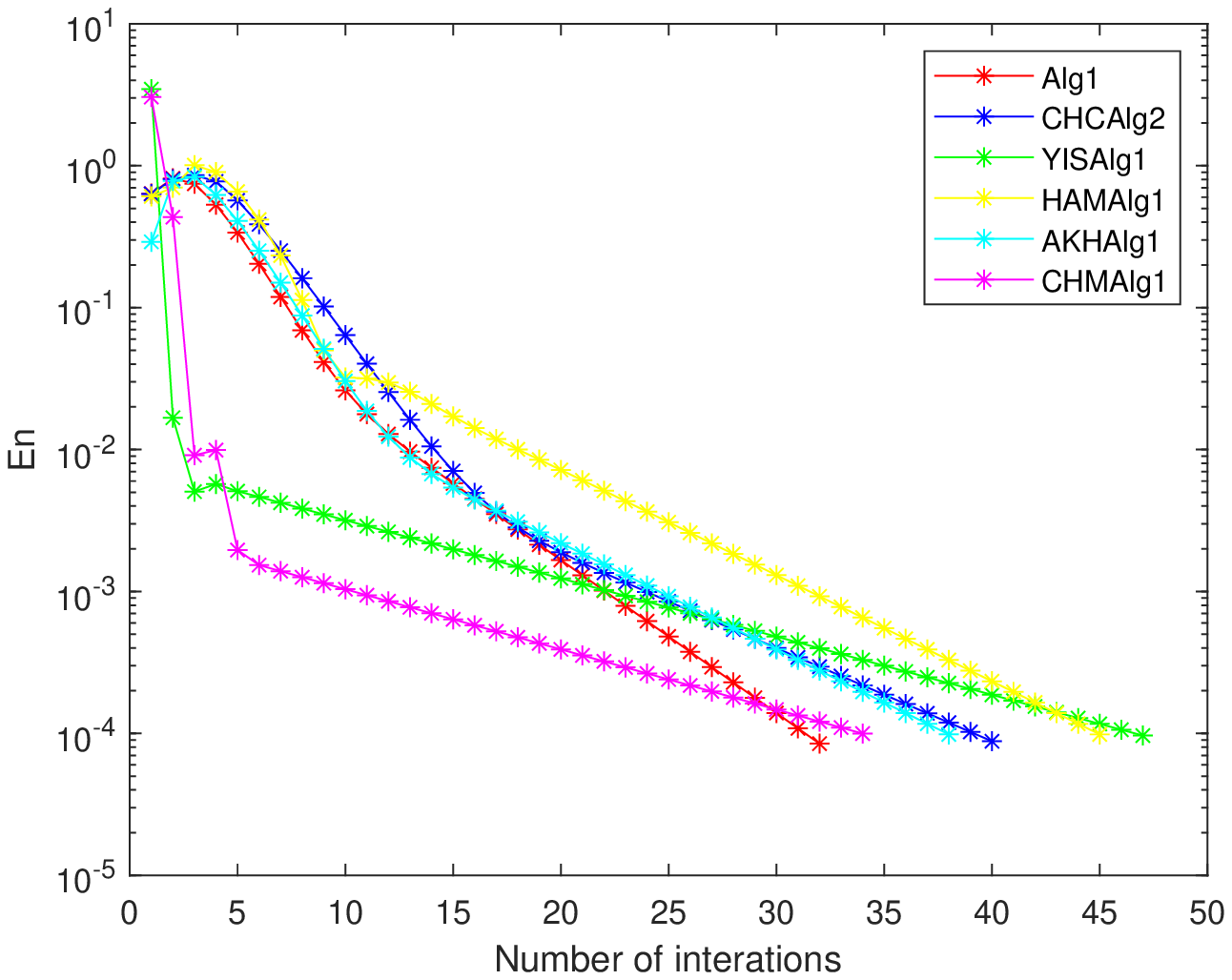}}
		\subfigure[${x_0} = \frac{{97{t^2} + 4t}}{{13}},{x_1} = \frac{{\sin (3t) + \cos (10t)}}{{100}}.$]{\includegraphics[width=7cm]{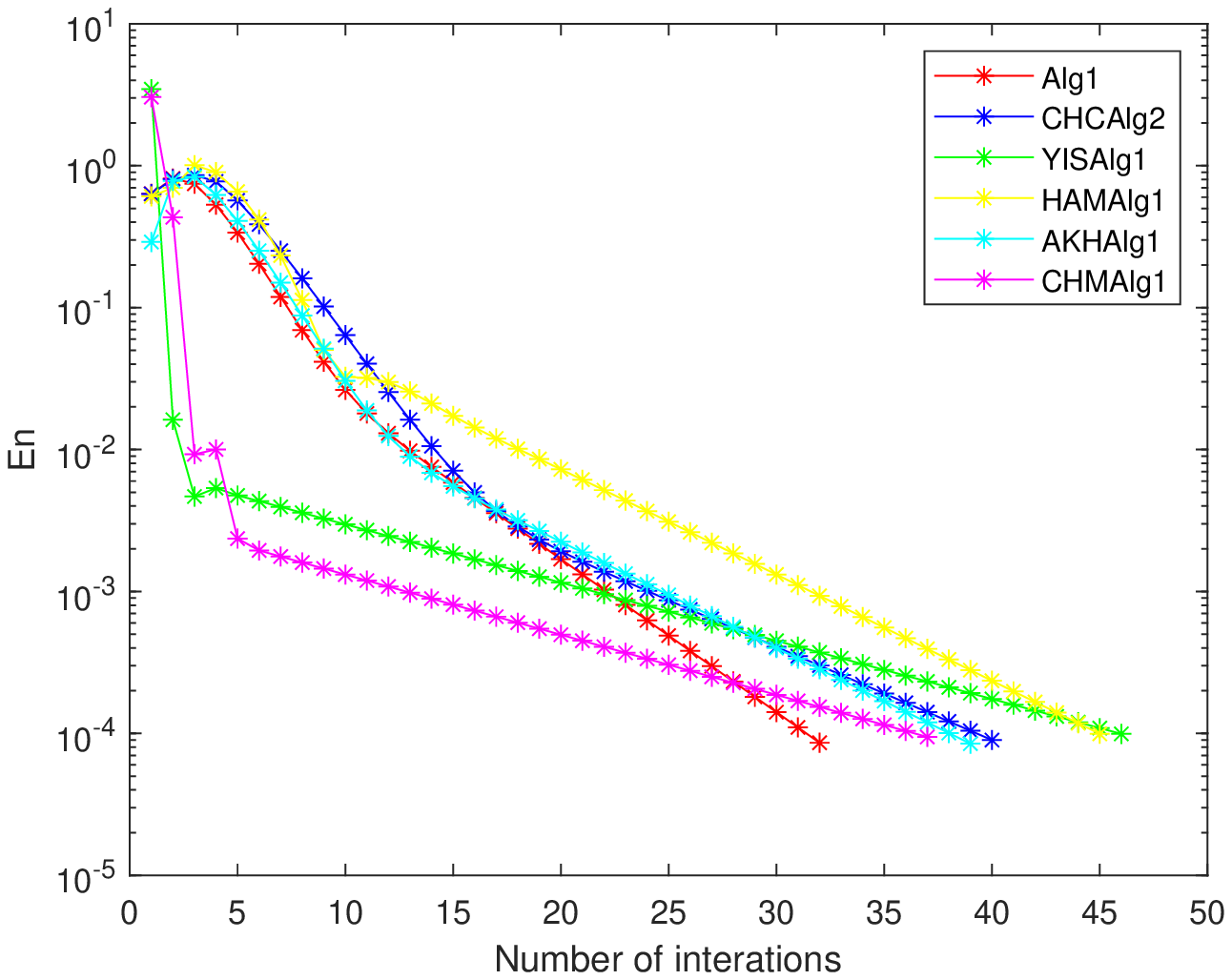}}
		\\ 
		\centering
		\subfigure[${x_0} = \frac{{{t^2} - {e^{ - 7t}}}}{{250}},{x_1} = \frac{{\sin (3t) + \cos (10t)}}{{100}}.$]{\includegraphics[width=7cm]{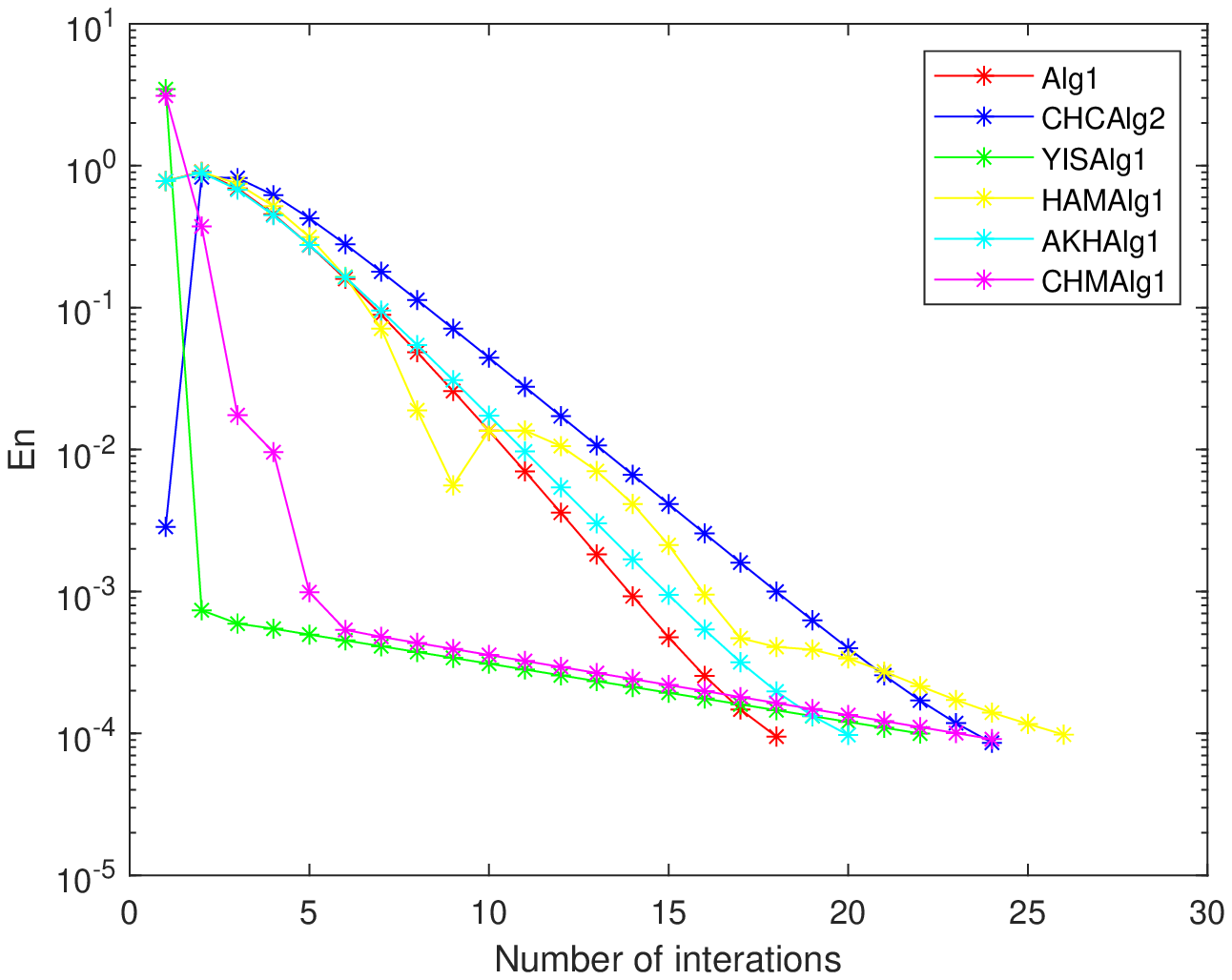}}
		\subfigure[${x_0} = \frac{{\sin (3t) + \cos (10t)}}{{100}},{x_1} = \frac{{97{t^2} + 4t}}{{13}}.$]{\includegraphics[width=7cm]{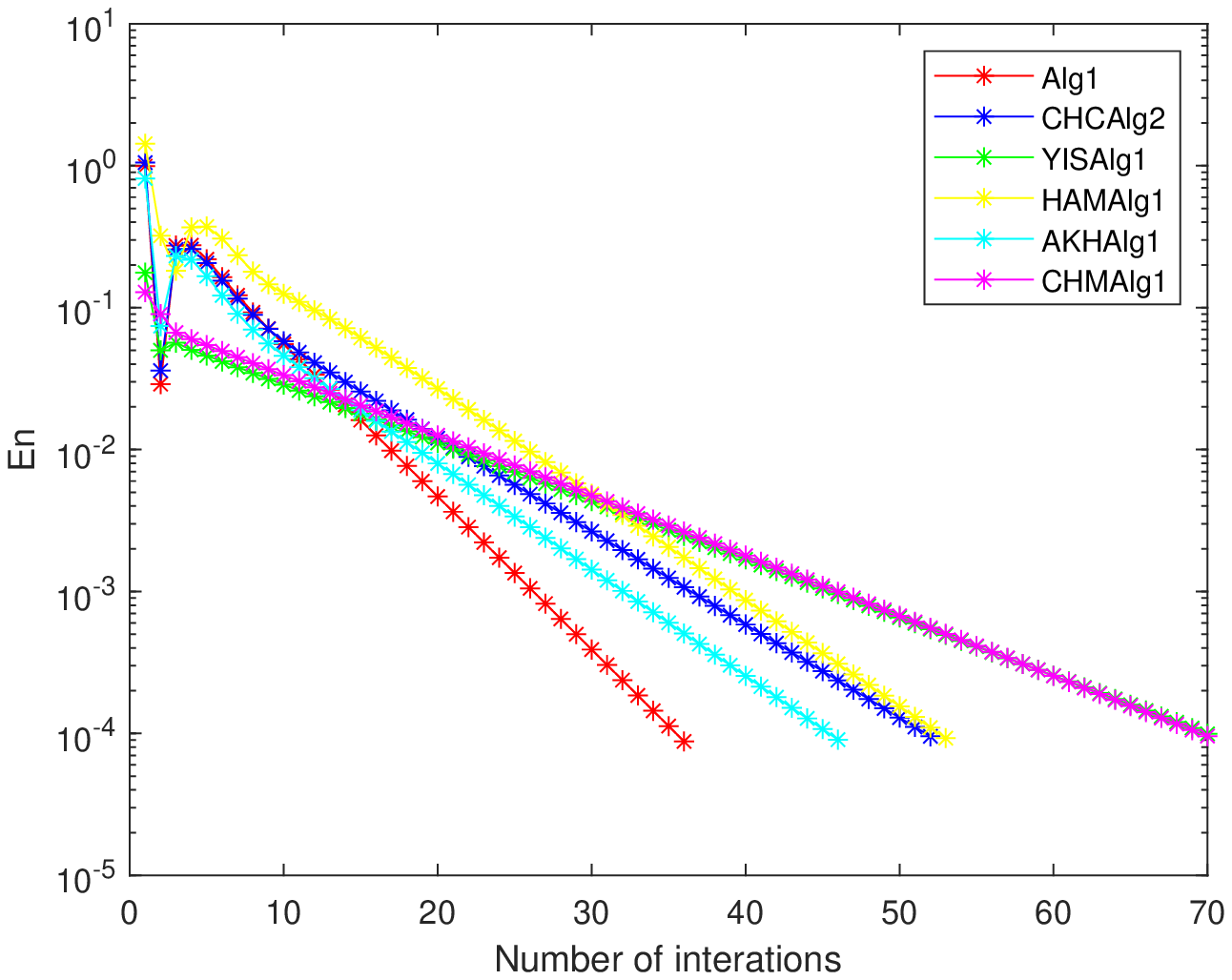}}
		\caption{ Numerical behavior of $E_n$ for Example \ref{E5.3}} 
	\end{figure}
	\begin{remark}
		\label{r5.3} From Table 5 and Fig.4, we observe that our Algorithm 3.1 performs better and converges
		faster.
	\end{remark}	
\end{example}
\section{Conclusions}

In this paper, we propose a new Tseng splitting method with double inertial extrapolation steps for solving monotone inclusion problems in real Hilbert spaces and establish the weak convergence, nonasymptotic $O(\frac{1}{\sqrt{n}})$ convergence rate, strong convergence and linear convergence rate of the proposed algorithm, respectively. Our method has the following advantages:
\par (i) Our method uses adaptive step sizes, which can be updated by a simple calculation without knowing the Lipschitz constant of the underlying operator. \par (ii) Our method own double inertial extrapolation steps, in which inertial factor $\alpha_n$ can equal $1$. This is not allowed in the corresponding algorithms of \cite{CV,CH}, where the only single inertial extrapolation step is considered and the inertial factor is bounded away from $1$. From Table 1 in section 5, it can be seen Algorithm 3.1 with double inertial extrapolation steps
 outperforms the one with the single inertia.
 \par (iii) Our method includes the corresponding methods considered in \cite{AB,CV,YI} as special cases. Especially, when our algorithm is used to solve variational inequalities, the relaxed parameter sequence $\{\theta_n\}$ have larger choosing interval than the ones of \cite{YI}. Via Table 3 in section 5, we observe that the performance of the algorithm becomes better with the increase of the relaxed parameter $\theta_n$.
 \par (iv) To the best of our knowledge, there are few available convergence rate results for algorithms with the double inertial extrapolation steps for solving variational inequalities and monotone inclusions. From numerical experiments in section 5, we can see that our algorithm has better efficiency than the corresponding algorithms in \cite{AB,CV,CH,VA,YI}.
\begin{acknowledgements}
This work was supported by the National Natural Science Foundation of China (11701479, 11701478), the Chinese Postdoctoral Science Foundation (2018M643434) and the Fundamental Research Funds for the Central Universities (2682021ZTPY040).
\end{acknowledgements}

\end{document}